\title{ ~~\\ The hexagonal versus the square lattice}
\author{Pieter Moree and Herman J.J. te Riele}
\documentclass[12pt]{article}
\usepackage{amssymb, latexsym, amsfonts, graphicx}
\def\@ptsize{2}
\newtheorem{Thm}{Theorem}
\newtheorem{Conj}{Conjecture}
\newtheorem{Lem}{Lemma}
\newtheorem{cor}{Corollary}

\newtheorem{Prop}{Proposition}
\newcommand{\qed}{\hfill $\Box$}

\begin{document}
\date{}
\maketitle
{\def\thefootnote{}
\footnote{\noindent P. Moree: KdV Institute, University of Amsterdam,
Plantage Muidergracht 24, 1018 TV Amsterdam, The Netherlands, e-mail:
moree@science.uva.nl\\
\indent H.J.J. te Riele: CWI, Postbus 94079, 1090 
GB Amsterdam, The Netherlands,\\ e-mail: herman@cwi.nl}}
{\def\thefootnote{}
\footnote{{\it Mathematics Subject Classification (2000)}.
11N13, 11Y35, 11Y60}}
\begin{abstract}
\noindent Schmutz Schaller's \cite[p. 201]{conjecture} conjecture
regarding the lengths of the hexagonal versus the lengths of the
square lattice is shown to be true.
The proof uses results from (computational) prime number theory and
from \cite{moree}.\\
\indent Using an identity due to Selberg, it is shown that the conjecture
can in principle be resolved also without using computational prime number
theory. By our approach, however, this would require a 
huge amount of computation,
\end{abstract}
\section{Introduction}
In \cite[p. 201]{conjecture} Schmutz Schaller, motivated
by considerations from hyperbolic geometry, makes the conjecture
that in dimensions 2 to 8 the best known lattice sphere packings have
`maximal lengths' and goes on to write: "In dimension 2 the conjecture
means in particular that the hexagonal
lattice is `better' than the square lattice. More precisely, let
$0<h_1<h_2<\cdots$ be the positive integers, listed in ascending order,
which can be written as $h_i=x^2+3y^2$ for integers $x$ and $y$.
Let $0<q_1<q_2<\cdots$ be the positive integers, listed in ascending
order, which can be written as $q_i=x^2+y^2$ for integers $x$ and $y$.
Then
the conjecture is that $q_i\le h_i$ for $i=1,2,3,\cdots$."
For an introduction to these matters see \cite{student}.
For some progress regarding Schmutz Schaller's general conjecture in
dimension 2 see \cite{kuhnlein} (this case of the conjecture
is also mentioned in \cite[p. xxx]{latticebible}).\\
\indent For $j\ge 1$ let $b_j(n)=1$ if $n$ is represented by the quadratic
form $X^2+jY^2$ and $b_j(n)=0$ otherwise. The characteristic 
functions $b_1$ and $b_3$ are well understood.
The following result was already known to Fermat.
\begin{Lem}
\label{classical}
A positive integer $n$ is represented by the form $X^2+Y^2$ if and only if every 
prime
factor $p$ of $n$ of the form $p\equiv 3({\rm mod~}4)$ occurs to an
even power. A positive integer $n$ is represented by the form
$X^2+3Y^2$ if and only if
every prime factor $p$ of $n$ of the form $p\equiv 2({\rm mod~}3)$ occurs
to an even power.
\end{Lem}
(In general the natural numbers $n$ that are represented by a quadratic form
$X^2+mY^2$ are rather more difficult to describe, cf. the beautiful
book of D. Cox \cite{cox}.) Lemma \ref{classical} implies
that $b_1$ and $b_3$ are multiplicative functions.\\
\indent Let $B_i(x)=\sum_{n\le x}b_i(n)$
for $i=1$ and $i=3$.
Schmutz Schaller's conjecture regarding
the square versus the hexagonal lattice can
be reformulated as follows in terms of $B_1$
and $B_3$.
\begin{Conj}
We have $B_1(x)\ge B_3(x)$ for every $x$.
\end{Conj}
\indent The first asymptotic
result on $B_1(x)$ goes back
to Landau \cite{landau},
who proved in 1908 that
\begin{equation}
\label{1908}
B_1(x)\sim C_{b_1}{x\over \sqrt{\log x}},
\end{equation}
where
\begin{equation}
\label{landaucons}
C_{b_1}={1\over \sqrt{2}}\prod_{p\equiv
3({\rm mod~}4)}(1-p^{-2})^{-1/2}
={\pi\over 4}\prod_{p\equiv 1({\rm mod~}4)}
(1-p^{-2})^{1/2}\approx 0.764.
\end{equation}
(Here and in the sequel the letter $p$ is used to indicate primes.)
Landau's proof uses contour integration. It is not difficult to use
his method to show, cf. \cite{serre}, that for every $k\ge 2$ there exists 
constants
$C_{b_1}(2),\cdots,C_{b_1}(k)$ such that
\begin{equation}
\label{poincare}
B_1(x)=C_{b_1}{x\over \sqrt{\log x}}\left(1+{C_{b_1}(2)\over \log x}
+\cdots+{C_{b_1}(k)\over \log^{k-1} x}+O\left({1\over \log^{k}x}\right)\right).
\end{equation}
This result can also be established by methods not using complex
analysis, cf. \cite[p. 288]{postnikov}.
At the beginning of 1913
a then unknown Hindu clerk by the
name of Ramanujan wrote in his first letter
to Hardy \cite{berndt} that he could prove that
\begin{equation}
\label{eerstebrief}
B_1(x)=C_{b_1}\int_2^x{dt\over \sqrt{\log t}}+O(x^{1-\varepsilon}),
\end{equation}
for some $\varepsilon>0$. (For a reconstruction of Ramanujan's speculative
argument see \cite[pp. 60-66]{berndtr}.)
Note the similarity of Ramanujan's claim with
the prime number theorem. From (\ref{eerstebrief}) we infer that
$C_{b_1}(2)=1/2$ by partial integration.
Shanks \cite{shanks} showed, however, that $C_{b_1}(2)\ne 1/2$, thus disproving 
Ramanujan's claim. (Ramanujan gave the correct formula and numerical
approximation for $C_{b_1}$ though.)
The constants $C_{b_1}$ and $C_{b_1}(2)$
are known as the {\it Landau-Ramanujan constant} and the
{\it second order
Landau-Ramanujan constant}, respectively. For more on the
evaluation of these constants see Section \ref{constants}. For
more on mathematical constants in general see e.g. \cite{finch}.\\
\indent Ramanujan \cite{bono} 
stated several claims similar to
(\ref{eerstebrief}) in his `unpublished' manu- script on the partition 
and
tau functions, see Section 6. All of them are disproved
in \cite{moreerama}. It can be shown, however, that in each case
Ramanujan's claims give the correct asymptotic main term.\\
\indent Similarly to (\ref{poincare}) it can be shown
that for
arbitrary $k\ge 2$ there exist constants $C_{b_3}(2),\cdots,C_{b_3}(k)$  such 
that
\begin{equation}
\label{poincare2}
B_3(x)=C_{b_3}{x\over \sqrt{\log x}}\left(1+{C_{b_3}(2)\over \log x}
+\cdots+{C_{b_3}(k)\over \log^{k-1}x}+O\left({1\over \log^{k}x}\right)\right),
\end{equation}
where
$$C_{b_3}={1\over \sqrt{2}}{1\over 3^{1/4}}\prod_{p\equiv
2({\rm mod~}3)}(1-p^{-2})^{-1/2}
={\pi 3^{1/4}\sqrt{2}\over 9}\prod_{p\equiv 1({\rm mod~}3)}
(1-p^{-2})^{1/2}\approx 0.639.$$
We thus arrive at the following conclusion.
\begin{Prop}
Conjecture {\rm 1} is asymptotically true.
\end{Prop}
\indent Table 1 (copied from 
\cite{shanksschmid} and verified by the second author) suggests that
Conjecture 1 is true for small $x$ as well. The literature thus provides us with
good indications that Conjecture 1 is true. The purpose of this paper
is going beyond this and to prove that Conjecture 1 is indeed true. 
\begin{table}
\begin{center}
\begin{tabular}{|c|c|c||c|c|c|}\hline
$x$&$B_1(x)$&$B_3(x)$&$x$&$B_1(x)$&$B_3(x)$\\ \hline\hline
$2^{1}$&2&1&$2^{14}$&4357&3645\\ \hline
$2^2$&3&3&$2^{15}$&8363&6993\\ \hline
$2^3$&5&4&$2^{16}$&16096&13456\\ \hline
$2^4$&9&8&$2^{17}$&31064&25978\\ \hline
$2^5$&16&14&$2^{18}$&60108&50248\\ \hline
$2^6$&29&25&$2^{19}$&116555&97446\\ \hline
$2^7$&54&45&$2^{20}$&226419&189291\\ \hline
$2^8$&97&82&$2^{21}$&440616&368338\\ \hline
$2^{9}$&180&151&$2^{22}$&858696&717804\\ \hline
$2^{10}$&337&282&$2^{23}$&1675603&1400699\\ \hline
$2^{11}$&633&531&$2^{24}$&3273643&2736534\\ \hline
$2^{12}$&1197&1003&$2^{25}$&6402706&5352182\\ \hline
$2^{13}$&2280&1907&$2^{26}$&12534812&10478044\\ \hline
\end{tabular}
\caption{$B_1(x)$ versus $B_3(x)$}
\end{center}
\end{table}
\begin{Thm}
\label{main}
We have $B_1(x)\ge B_3(x)$ for every $x$, that is 
Schmutz Schaller's conjecture that the hexagonal lattice is `better'
than the square lattice is true.
\end{Thm}
\indent Landau's classical result (\ref{1908}) has been generalised in many
directions, see \cite{cazaran} for a survey with over 50 references.
Despite this rich history, nobody but
the present author (in \cite{moree}) seems to have been
concerned with proving effective results in this area, which is precisely
what is needed to establish Theorem \ref{main}.\\
\section{Preliminaries}
Let $f$ be a multiplicative function from the natural numbers to $\mathbb R_{\ge 
0}$.
We define $M_f(x)=\sum_{n\le x}f(n)$,
$\mu_f(x)=\sum_{n\le x}f(n)/n$
and $\lambda_f(x)=\sum_{n\le x}f(n)\log n$.
We denote the formal Dirichlet series $\sum_{n=1}^{\infty}f(n)n^{-s}$
associated to $f$ by $L_f(s)$.
We define $\Lambda_f(n)$
by
$$-{L_f'(s)\over L_f(s)}=\sum_{n=1}^{\infty}{\Lambda_f(n)\over n^s}.$$
Notice that
\begin{equation}
\label{convolutie}
f(n)\log n=\sum_{d|n}f(d)\Lambda_f({n\over d}).
\end{equation}
If $f$ is the characteristic function of a subsemigroup of the natural
integers with $(1<)q_1<q_2<\cdots$ as generators, then it can be shown
that $\Lambda_f(n)=\log q_i$ if $n$ equals a positive power of a generator
$q_i$ and $\Lambda_f(n)=0$ otherwise.
For $f=b_1$ we thus find, using Lemma \ref{classical},
$$
\Lambda_{b_1}(n)=\cases{
2\log p &if $n=p^{2r},~r\ge 1$ and $p\equiv 3({\rm mod~}4)$;\cr
\log p &if $n=p^r,~r\ge 1$ and $p\equiv 1({\rm mod~}4)$ or $p=2$;\cr
0 &otherwise.}
$$
For $f=b_3$ we find
$$
\Lambda_{b_3}(n)=\cases{
2\log p &if $n=p^{2r},~r\ge 1$ and $p\equiv 2({\rm mod~}3)$;\cr
\log p &if $n=p^r,~r\ge 1$ and $p\equiv 1({\rm mod~}3)$ or $p=3$;\cr
0 &otherwise.}
$$
From property (\ref{convolutie}) of
$\Lambda_f(n)$, we easily infer that
\begin{equation}
\label{lnaarpsi}
\lambda_f(x)=\sum_{n\le x}f(n)\psi_f({x\over n}),
\end{equation}
where $\psi_f(x)=\sum_{n\le x}\Lambda_f(n)$. The function $\Lambda_f$
and $\psi_f$ are analogues of, respectively, the von Mangoldt and the
Chebyshev $\psi$-function.

\section{Some related conjectures}
Unfortunately it seems that $M_f$ is not a very
`natural' mathematical object, $\mu_{f}$ is on
the other hand (as is amply demonstrated by browsing
through the literature). For this reason we
consider two additional conjectures:
\begin{Conj}
\label{conjecture2}
We have $\lambda_{b_1}(x)\ge \lambda_{b_3}(x)$ for $x\ge 8$.
\end{Conj}
\begin{Conj}
\label{conjecture3}
We have $\mu_{b_1}(x)\ge \mu_{b_3}(x)$ for every $x$.
\end{Conj}
Note that $\exp(\lambda_{b_1}(x)/2)$ is the product of all different
lengths in the square lattice not exceeding $\sqrt{x}$. Thus Conjecture 2
can be reformulated as stating that the product of the different
distances not exceeding $x$ occuring in the square lattice always exceeds the 
product of the
different distances not exceeding $x$ in the hexagonal lattice, provided
that $x\ge 2\sqrt{2}$.\\
\indent Conjecture 1 clearly implies Conjecture 3.
Furthermore we have:
\begin{Prop}
\label{proposition2}
Conjecture {\rm 2} implies Conjecture {\rm 1}.
\end{Prop}
{\it Proof}. We have, for $x\ge 2$,
\begin{equation}
\label{mlambda}
M_f(x)=\int_{2-}^x{d\lambda_f(t)\over \log t}=
{\lambda_f(x)\over \log x}+\int_2^x {\lambda_f(t)\over t\log^2 t}dt.
\end{equation}
Denote the latter integral by $I_f(x)$. It is not difficult to
show that $I_{b_1}(x)\ge I_{b_3}(x)$ for $x\le 8$. Conjecture 2
then implies that the latter inequality holds for every $x$. The
truth of
Conjecture 2 together with (\ref{mlambda}) then implies
that $B_1(x)\ge B_3(x)$ for $x\ge 8$. By direct computation
we then infer that the latter inequality holds for every $x$. \qed

Thus in order to establish Theorem \ref{main}, it suffices to
establish Conjecture 2. 
From (\ref{lnaarpsi}) and $\psi_{b_i}(x)\sim x/2$ as $x$ tends to
infinity it follows that $\lambda_{b_i}(x)\sim \mu_{b_i}(x)/2$, as
$x$ tends to infinity. An effective form of this relationship, together
with an effective estimate for $\mu_{b_i}$ (provided
by Lemma \ref{lemma1}), then allows us to prove the
main result of this paper:
\begin{Thm}
\label{hiergaathetom}
The Conjectures {\rm 1}, {\rm 2} and {\rm 3} are all true.
\end{Thm}
Complications arise due to the fact that
$\lim_{x\rightarrow \infty}B_1(x)/B_3(x)=1.1961377420\cdots$, which is rather 
close to 1 
and
that $\psi_{b_i}(y)$ is not so close to $y/2$ for various 
ranges of small $y$ (the convolutional nature of (\ref{lnaarpsi}) forces
us to take the small $y$ range into account).

\section{The toolbox}
The following result from \cite{moree} will play a crucial r\^ole. It is in 
essence an effective version of Theorem A of \cite{song}. 
\begin{Lem}
\label{lemma1}
Let $f$ be a multiplicative function from the natural numbers to $\mathbb R_{\ge 
0}$.
Suppose that there exist constants $D_{-}$,
$D_{+}$ and $\tau$, with $\tau>0$, such that for every $x\ge x_0$,
\begin{equation}
\label{volgendeinklemming}
D_{-}\mu_f(x)\le 
\sum_{n\le x}{f(n)\over n}\left\{\sum_{m\le {x\over n}}{\Lambda_f(m)\over m
}-\tau
\log{x\over n}\right\}
\le D_{+}\mu_f(x).
\end{equation}
Then we have, for $x>\max\{x_0,\exp(D_{+})\}$,
\begin{equation}
\label{secondestimate}
{C_f\over \tau}\log^{\tau}x{\left(1-{D_+\over \log x}\right)^{\tau+1}\over
1-{D_{-}\over \log x}}
\le \mu_f(x)
\le {C_f\over \tau}\log^{\tau}x
{\left(1-{D_{-}\over \log x}\right)^{\tau+1}\over
1-{D_{+}\over \log x}},
\end{equation}
where
\begin{equation}
\label{defcee}
C_f:={1\over \Gamma(\tau)}\lim_{s\rightarrow 1+0} (s-1)^{\tau}L_f(s).
\end{equation}
In particular, if there exist constants
$C_{-}$ and $C_{+}$
such that
\begin{equation}
\label{inklemming}
C_{-}\le \sum_{n \le x}{\Lambda_f(n)\over n}-\tau \log
x \le C_{+}~~{\rm for~}x\ge 1,
\end{equation}
then {\rm (\ref{secondestimate})} holds true, for $x>\exp(C_{+})$,
with $D_{-}=C_{-}$ and $D_{+}=C_{+}$.
\end{Lem}
\noindent {\tt Remark 1}. From 
the proof of this lemma, $\int_1^x{\mu_f(t)dt/t}$
appears as a more easily estimated function than $\mu_f(x)$.
Interestingly, Landau \cite{landau} in his proof
of (\ref{1908}) using contour integration, estimates
$\int_1^x{\mu_{b_1}(t)dt/t}$ rather
than $B_1(x)$ itself.\\
{\tt Remark 2}. If $\lim_{x\rightarrow \infty}(\sum_{n\le x}(\Lambda_f(n)/n
-\tau \log x))$ exists, we denote this by $B_f$.\\

\noindent Let us put 
$$L(x,\tau,D_{-},D_{+})={(\log x -D_{+})^{\tau+1}\over \log x -D_{-}}
{\rm ~and~}
U(x,\tau,D_{-},D_{+})={(\log x-D_{-})^{\tau+1}\over \log x -D_{+}}.$$ 
Thus we can write (\ref{secondestimate}) as
$C_fL(x,\tau,D_{-},D_{+})/\tau\le 
\mu_f(x)\le C_fU(x,\tau,D_{-},D_{+})/\tau$.\\
\indent Let $r,s$ and $c_1$ be given. At a few instances in the sequel
we want to show that for  
every $x\ge x_2$, with $x_2$ some explicit constant,
we have
$\mu_f(x/r)\ge c_1\mu_g(x/s)$, where $g$ satisfies the
conditions of Lemma \ref{lemma1} with constants $\tau$, $D'_{-}$ and
$D'_{+}$.
By Lemma \ref{lemma1} this leads us to
consider inequalities of the form
\begin{equation}
\label{veelvariabelen}
L({x\over r},\tau, D_{-}, D_{+})\ge c_2U({x\over s},\tau, D'_{-}, D'_{+})
\end{equation}
where all variables and constants are real numbers with $\tau,r,s$ and
$c_2$ positive, $D_{-}\le D_{+}$, $D'_{-}\le D'_{+}$ and
$x\ge x_0:=\max\{\exp(D'_{+})s,\exp(D_{+})r\}$.
We recall the following lemma from \cite{moree}:
\begin{Lem}
\label{delaatstehoopik}
If $\log s+D'_{-}\le D_{+}+\log r$ and {\rm (\ref{veelvariabelen})} is satisfied
for some $x_1>x_0$, then {\rm (\ref{veelvariabelen})} is satisfied for every
$x\ge x_1$. If $\log s+D'_{-}>D_{+}+\log r$, and 
$$c_2\left(1+{D'_{+}-D'_{-}\over \log(x_1/s)-D'_{+}}\right)
\le 1+{D_{-}-D_{+}\over \log(x_1/r)-D_{-}},$$
for some $x_1>x_0$, then
{\rm (\ref{veelvariabelen})} is satisfied for every $x\ge x_1$.
\end{Lem} 
We also need the following result about the difference
between $U({x\over r},{1\over 2}, D_{-}, D_{+})$ 
and $L({x\over s},{1\over 2}, D_{-}, D_{+})$.
\begin{Lem}
\label{dalertje}
Assume that $D_{+}>D_{-}$ and $s\ge r\ge 1$.
The difference 
$$U({x\over r},{1\over 2},D_{-},D_{+})-L({x\over s},{1\over 2},D_{-},D_{+})$$ 
is
monotonically decreasing for $x\ge s\exp(1.01D_{+}-0.01D_{-})$.
\end{Lem}
The difference in the latter lemma multiplied by $C_{b_i}$
appears if we try to bound above $\mu_{b_i}(x/r)-\mu_{b_i}(x/s)$. 
Notice that the latter difference is not monotonically
decreasing from
any $x$ onwards, although
it can be bounded above by a function that is monotonically decreasing for
all sufficiently large $x$.\\
\indent Our proof of Lemma \ref{dalertje} uses the following lemma. 
\begin{Lem}
\label{beetjegek}
Let $y$ and $\delta$ be non-negative real numbers. Then
the inequality
\begin{equation}
\label{welraar}
\sqrt{y+1+\delta}(y+\delta-2)(y+1)^2\le \sqrt{y}(y+3)(y+\delta)^2,
\end{equation}
holds if either $\delta\le 2$ or $y\ge 0.0099945$.
\end{Lem}
{\it Proof}. On replacing the inequality sign in
(\ref{welraar}) with the equality sign and
squaring both sides, we obtain an equation of an
algebraic curve. Using continuity
and e.g. Maple's function {\tt fsolve} (for numerically determining roots
of polynomial equations), the result can then be deduced. \qed\\

\noindent {\tt Remark}. For $y=0.0099944$ and $\delta\approx 5.4$ inequality 
(\ref{welraar})
is not satisfied. Indeed if we square both sides of the inequality and
take the difference then, considered as a polynomial in $y$, the
discriminant has $27\delta^5-198\delta^4+410\delta^3-936\delta^2+1299\delta-730$ 
as a factor, which has $5.44694735\cdots$ as 
its largest real root.
Considered as a polynomial in $d$, we find
$$27y^8-72y^7-2380y^6-12792y^5-33822y^4-48888y^3-32076y^2-2376y+27$$ as
a factor of the discriminant, which has $0.00999445028\cdots$ as its one but 
second largest
real root.\\

\indent We can now prove Lemma \ref{dalertje}.\\

\noindent {\it Proof of Lemma} \ref{dalertje}. Differentiating
$U(x/r,{1\over 2},D_{-},D_{+})-
L(x/s,{1\over 2},D_{-},D_{+})$ yields, after some tedious calculations,
that the derivative is non-positive provided that
(\ref{welraar}) is satisfied with $y=(\log(x/s)-D_{+})/(D_{+}-D_{-})$
and $\delta=\log(s/r)/(D_{+}-D_{-})$. The result then follows on
invoking Lemma \ref{beetjegek}. \qed\\

\section{Numerical evaluation of certain constants}
\label{constants}
For our proof of Theorem \ref{main} we need to evaluate
the constants $C_{b_1},C_{b_3},B_{b_1}$ and $B_{b_3}$ with
enough numerical precision. The purpose of this section is
to achieve this. (Recall that
$B_{b_i}=\lim_{x\rightarrow\infty}(\sum_{n\le x}\Lambda_{b_i}(n)/n-(\log 
x)/2)$.)\\
\indent We first consider the evaluation of $C_{b_3}$ and
$C_{b_1}$ (defined by (\ref{defcee})). We have, for Re$(s)>1$,
$$L_{b_3}(s)=(1-3^{-s})^{-1}\prod_{p\equiv 1({\rm mod~}3)}
(1-p^{-s})^{-1}\prod_{p\equiv 2({\rm mod~}3)}
(1-p^{-2s})^{-1}$$
and
\begin{equation}
\label{generating}
L_{b_3}(s)^{2}=\zeta(s)L(s,\chi_{-3})(1-3^{-s})^{-1}
\prod_{p\equiv 2({\rm mod~}3)}
(1-p^{-2s})^{-1}.
\end{equation}
From this, (\ref{defcee}), $\lim_{s\rightarrow 1+0}(s-1)\zeta(s)=1$ and
the fact that $\Gamma({1\over 2})=\sqrt{\pi}$, we obtain
$$C_{b_3}^2={3L(1,\chi_{-3})\over 2\pi}
\prod_{p\equiv 2({\rm mod~}3)}
(1-p^{-2})^{-1},$$
where for any fundamental discriminant $D$, $\chi_D$ denotes 
Kronecker's extension $(D/n)$ of
the Legendre symbol \cite[Chapter 5]{davenport}.
If $\chi$ is a real primitive character modulo $k$ and $\chi(-1)=-1$,
then
$$L(1,\chi)=-{\pi\over k^{3/2}}\sum_{n=1}^k n\chi(n),$$
by Dirichlet's celebrated class number
 formula (cf.
equation (17) of \cite[Chapter 6]{davenport}).
We infer that $L(1,\chi_{-3})=\pi/\sqrt{27}$.
Using that $C_{b_3}$ must be positive and $\zeta(2)=\pi^2/6$, we
then infer that
$$C_{b_3}={1\over \sqrt{2}}{1\over 3^{1/4}}\prod_{p\equiv
2({\rm mod~}3)}(1-p^{-2})^{-1/2}
={\pi 3^{1/4}\sqrt{2}\over 9}\prod_{p\equiv 1({\rm mod~}3)}
(1-p^{-2})^{1/2}.$$
Likewise, using that $L(1,\chi_{-4})=\pi/4$, we find formula
(\ref{landaucons}) for $C_{b_1}$.\\
\indent Note that, for $\Re (s)>1/2$,
\begin{equation}
\label{shanksje}
\prod_{p\equiv 3({\rm mod~}4)}(1-p^{-2s})^{-2}
={\zeta(2s)(1-2^{-2s})\over L(2s,\chi_{-4})}\prod_{p\equiv 3({\rm mod~}4)}
(1-p^{-4s})^{-1}.
\end{equation}
By recursion we then find from (\ref{landaucons}) and
(\ref{shanksje}) the following formula, 
$$C_{b_1}={1\over \sqrt{2}}\prod_{n=1}^{\infty}\left((1-2^{-2^n})
{\zeta(2^n)\over L(2^n,\chi_{-4})}\right)^{1/2^{n+1}},$$
which was already known to  Ramanujan
\cite[pp. 60-66]{berndtr} and
by Shanks \cite[p. 78]{shanks}.
Using this expression one 
computes that
$C_{b_1}=0.76422365358922066299\cdots$
Similarly one can show that
$$C_{b_3}={1\over \sqrt{2}}{1\over 3^{1/4}}\prod_{n=1}^{\infty}
\left((1-3^{-2^n})
{\zeta(2^n)\over L(2^n,\chi_{-3})}\right)^{1/2^{n+1}},$$
and use it to compute that
$C_{b_3}=0.63890940544534388225\cdots$
which is in agreement with the
first seven (out of eight) decimals computed for $C_{b_3}$
by Shanks and Schmid \cite{shanksschmid}.\\
\indent On noting that, for Re$(s)\ge 1$,
$$\sum_{n=1}^{\infty}{\Lambda(n)-1\over n^s}
=-{\zeta'(s)\over \zeta(s)}-\zeta(s),$$
and using that $\zeta(s)=1/(s-1)+\gamma+O(s-1)$, where $\gamma$
denotes Euler's constant, is the Taylor series for $\zeta(s)$ around
$s=1$, one infers that
\begin{equation}
\label{lambdaovern}
\sum_{n\le x}{\Lambda(n)\over n}=\sum_{n\le x}{1\over n}
-2\gamma+o(1)=\log x -\gamma+o(1).
\end{equation}
Taking the logarithmic derivative of
(\ref{generating}) one obtains that
$$-2{L_{b_3}'(s)\over L
_{b_3}(s)}=
-{\zeta'(s)\over \zeta(s)}-{L'(s,\chi_{-3})\over L(s,\chi_{-3})}
+{\log 3\over 3^s-1}+2\sum_{p\equiv 2({\rm mod~}3)}{\log p\over
p^{2s}-1},$$
from which one easily infers that
$$2\sum_{n\le x}{\Lambda_{b_3}(n)\over n}
=\sum_{n\le x}{\Lambda(n)\over n}-{L'(1,\chi_{-3})\over L(1,\chi_{-3})}
+{\log 3\over 2}+2\sum_{p\equiv 2({\rm mod~}3)}{\log p\over
p^{2}-1}+o(1),$$
which yields, on invoking (\ref{lambdaovern}),
$$2B_{b_3}=-\gamma-{L'(1,\chi_{-3})\over L(1,\chi_{-3})}+{\log 3\over 2}
+2\sum_{p\equiv 2({\rm mod~}3)}{\log p\over p^2-1}.$$
Similarly we deduce that
$$2B_{b_1}=-\gamma-{L'(1,\chi_{-4})\over L(1,\chi_{-4})}
+{\log 2}+2\sum_{p\equiv 3({\rm mod~}4)}{\log p\over p^2-1}.$$
\indent Note that the argument above yielded that
$$B_{b_i}=-\lim_{s\rightarrow 1+0}\left({L'_{b_i}(s)\over L_{b_i}(s)}
+{1\over 2(s-1)}\right).$$
This can be alternatively deduced from Serre's \cite{serre} proof
of (\ref{poincare}), cf. \cite{moreerama}.\\
\indent As to the numerical evaluation of, for
example, the latter prime sum we note that
$$2\sum_{p\equiv 3({\rm mod~}4)}{\log p\over p^2-1}=-
{d\over ds}\log \prod_{p\equiv 3({\rm mod~}4)}\left({1\over 1-p^{-2s}}
\right)\Big|_{s=1}.$$
Then, applying (\ref{shanksje}) $m$ times, we obtain
$$
\sum_{p\equiv 3({\rm mod~}4)}{\log p\over p^2-1}=
\sum_{p\equiv 3({\rm mod~}4)}{\log p\over p^{2^{m+1}}-1}
+{1\over 2}\sum_{n=1}^m\left\{{L'(2^m,\chi_{-4})\over L(2^m,\chi_{-4})}
-{\zeta'(2^m)\over \zeta(2^m)}-{\log 2\over 2^{2^m}-1}\right\}.$$
Similarly we have
$$
\sum_{p\equiv 2({\rm mod~}3)}{\log p\over p^2-1}=
\sum_{p\equiv 2({\rm mod~}3)}{\log p\over p^{2^{m+1}}-1}
+{1\over 2}\sum_{n=1}^m\left\{{L'(2^m,\chi_{-3})\over L(2^m,\chi_{-4})}
-{\zeta'(2^m)\over \zeta(2^m)}-{\log 3\over 3^{2^m}-1}\right\}.$$
Using these expressions, one computes that
$B_{b_1}=0.163897318634581595856\cdots$ and similarly
$B_{b_3}=0.1535522449949958272447\cdots$.\\
\indent Now we can invoke \cite[Theorem 4]{moree} to compute 
the constants $C_{b_1}(2)$ and $C_{b_3}(2)$.
They are 
given by $C_{f}(2)=(1+B_f)/2$ for $f\in \{b_1,b_3\}$.
We thus find
$$C_{b_1}(2)=0.581948659317290797928\cdots,~C_{b_3}(2)
=0.576776122497497913622\cdots$$
In \cite{shanksschmid}
the authors write (in my notation)
`$B_3(x)$ remains so closely proportional to $B_1(x)$ that it is not
clear from this data whether $C_{b_3}(2)>C_{b_1}(2)$ or
$C_{b_1}(2)<C_{b_3}(2)$. It would be
unlikely that they are exactly equal.' We thus have resolved this matter.\\
\indent The numerical data from Table 1 in conjunction with the
values of $C_{b_1}$, $C_{b_1}(2)$ and (\ref{poincare}) suggest that
$C_{b_1}(3)>0$ and $C_{b_1}(4)<0$. Similarly it seems plausible that
$C_{b_3}(3)>0$ and $C_{b_3}(4)<0$.\\

\section{Intermezzo: On a claim of Ramanujan}
In the previous section we have seen that $B_{b_3}<\log \sqrt{3}$. This
 knowledge suffices to disprove a claim that was made
in a celebrated, hitherto unpublished, manuscript of Ramanujan \cite{bono} on
the partition and tau-functions.\\
\indent Let $\tau$ denote Ramanujan's tau-function. Put $T_n=0$ if 
$3|\tau(n)$ and
$T_n=1$ otherwise. In Ramanujan's manuscript we read \cite[p. 64]{bono}:
`We can show by transcendental methods that
\begin{equation}
\label{sriniclaim}
\sum_{k=1}^n T_k={C\over 3}\int_1^n {dx\over \sqrt{\log x}}
+O\left({n\over (\log n)^r}\right),
\end{equation}
where $r$ is any positive number and
$$C={2^{1/2}\over 3^{1/4}}.{1-7^{-2}\over 1-7^{-3}}.{1-13^{-2}\over 1-13^{-3}}
.{1-19^{-2}\over 1-19^{-3}}\cdots {1\over \{(1-2^{-2})(1-5^{-2})(1-11^{-2})
\cdots \}^{1/2}},$$
$2,5,11,\cdots$ being primes of the form $3k-1$ and $7,13,19,\cdots$ being
primes of the form $3k+1$'. This implies that for almost all
$n$, $\tau(n)$ is divisible by 3.\\
\indent Using that $\tau(n)\equiv n\sigma_1(n)({\rm mod~}3)$, where 
$\sigma_1(n)$ denotes the sum of the positive divisors of $n$, it is
easy to see that $T_n$ is multiplicative and that 
\begin{equation}
\label{doorrama}
\sum_{n=1}^{\infty}{T(n)\over n^s}=\prod_{p\equiv 2({\rm  mod~}3)}
{1\over 1-p^{-2s}}\prod_{p\equiv 1({\rm mod~}3)}{1+p^{-s}\over 1-p^{-3s}}.
\end{equation}
From (\ref{doorrama}) it is not difficult to verify Ramanujan's claim
regarding the value of $C$.
By logarithmic differentiation we obtain from (\ref{doorrama}) that 
$$\sum_{n=1}^{\infty}{\Lambda_T(n)\over n^s}=
\sum_{p\equiv 2({\rm mod~}3)}{2\log p\over p^{2s}-1}
+\sum_{p\equiv 1({\rm mod~}3)}\left[{\log p\over p^s+1}+
{3\log p\over p^{3s}-1}\right].$$
On comparing this series with that for 
$\sum_{n=1}^{\infty}\Lambda_{b_3}(n)n^{-s}$,
it is easily seen, on using that $B_{b_3}<\log \sqrt{3}$, that
$$B_T=B_{b_3}-\sum_{p\equiv 1({\rm mod~}3)}{(2p+1)\log p\over
(p^2+p+1)(p+1)}-\log \sqrt{3}<B_{b_3}-\log \sqrt{3}<0,$$
indeed we have $B_T=-0.53\cdots$. This shows that
$$\sum_{k=1}^n T_k={C\over 3}{n\over \sqrt{\log n}}
\left(1+{0.23\cdots\over \log n}+O\left({1\over (\log 
n)^{1+\epsilon}}\right)\right),
$$
where $0.23\cdots=(1+B_T)/2\ne 0.5$ (here 
we invoked Theorem 4 of \cite{moree}) and $\epsilon>0$. Thus the above claim of 
Ramanujan is false for
every $r>3/2$ and true for $r\le 3/2$.\\
\indent Note that if it would be true that $B_T=0$, an amazing
identity for Euler's constant would result. The manuscript \cite{bono}
contains several further assertions of the type
(\ref{sriniclaim}) (with 3 replaced
by various other primes), all of which
are disproved for $r>3/2$ in \cite{moreerama}.

\section{On the behaviour of $\sum_{n\le x}{\Lambda_{b_i}(n)\over n}-{\log 
x\over 2}$}
 Put
$H_i(x)=\sum_{n\le x}\Lambda_{b_i}(n)/n-\log \sqrt{x}$, for $i=1$ and
$i=3$. A good understanding of the behavior of $H_i$ is needed
in order to apply 
our key lemma, Lemma \ref{lemma1}. Let us define for $i=1$ and $i=3$,
$C_{+}(b_i)=\sup_{x\ge 1}H_i(x)$ and $C_{-}(b_i)=\inf_{x\ge 1}H_i(x)$.
\\ \\
${~~~~~~~~}${\resizebox{5.0cm}{!}{\includegraphics{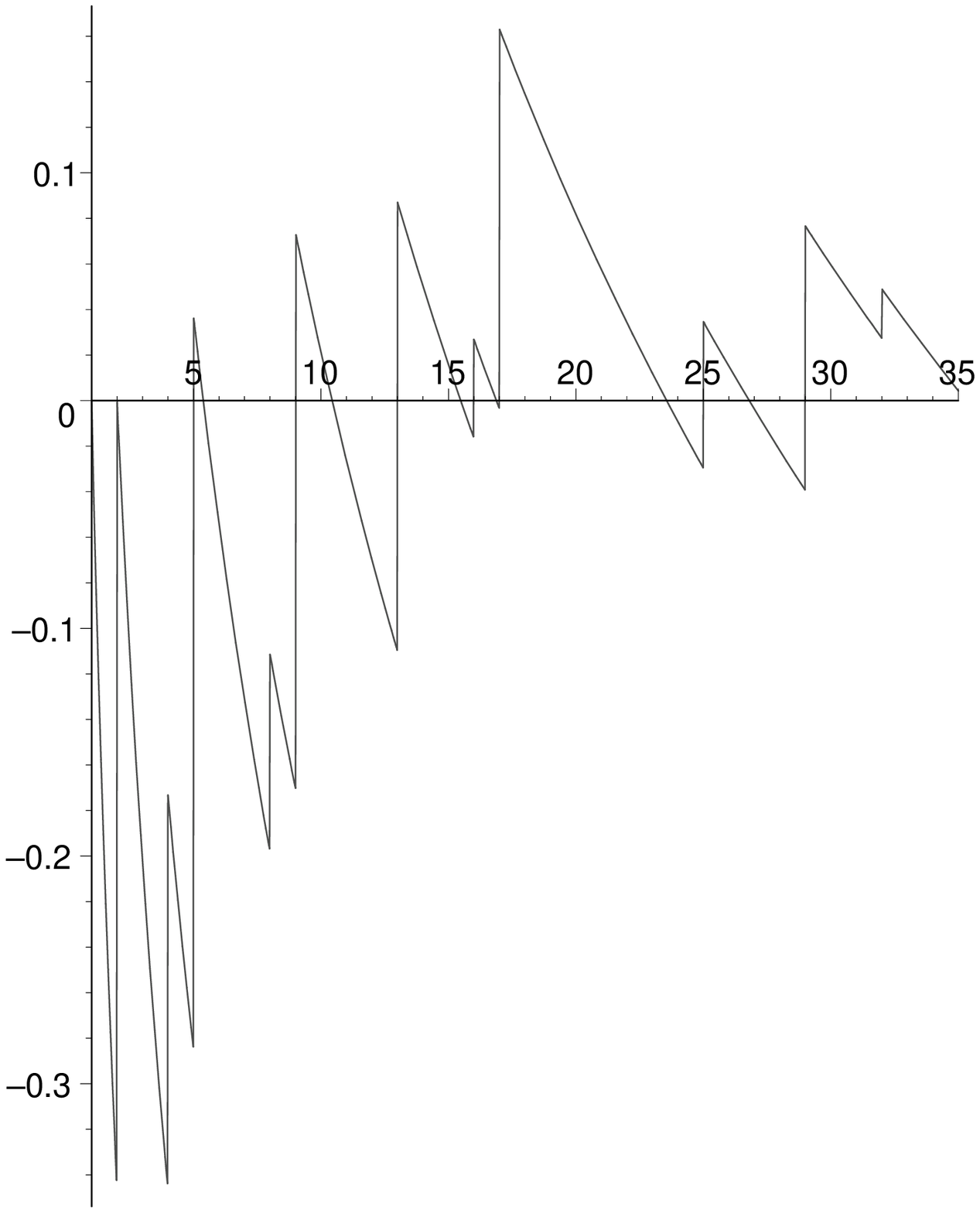}}}
${~~~~~~~~~~~~~}${\resizebox{5.0cm}{!}{\includegraphics{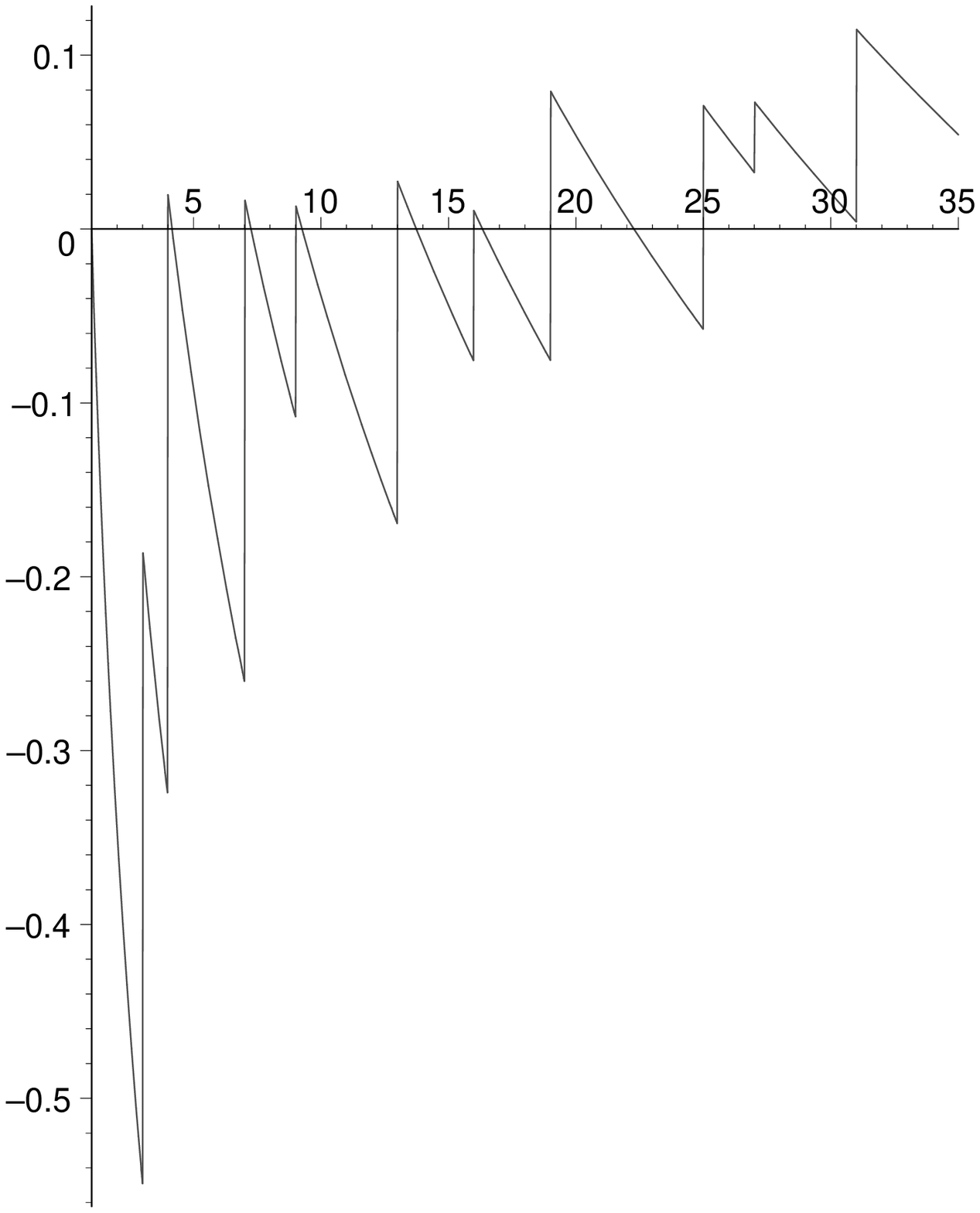}}}\\
${~~~~~~~~~~~~~~~~~~~~~~}${\bf Figure 1} 
${~~~~~~~~~~~~~~~~~~~~~~~~~~~~~~~~~~~~}$ {\bf Figure 2}\\
${~~~~~~~}${\tt Plot of $H_1(x)$ for $1\le x\le 35$}
${~~~~~~~}${\tt Plot of $H_3(x)$ for $1\le x\le 35$}
\\ \\
\noindent As in \cite{mertens} we define $V(x;d,a)=\sum_{n\le x\atop
n\equiv a({\rm mod~}d)}{\Lambda(n)/n}$. It can be shown that as $x$ tends
to infinity $V(x;d,a)-\log x/\varphi(d)$ tends to a limit $C(d,a)$.
Ramar\'e \cite{mertens} has established the following result.
\begin{Thm} 
\label{olivier}
{\rm \cite{mertens}}. For $x\ge 68$ we have
$|V(x;3,1)-{1\over 2}\log x-C(3,1)|\le 0.1205$ and 
$|V(x;4,1)-{1\over 2}\log x-C(4,1)|\le 0.0961$.
\end{Thm}
We recall from \cite{moree} that
\begin{equation}
\label{tochnodig}
\sum_{p^r>\sqrt{x}\atop p\equiv a({\rm mod~}d)}{\log p\over p^{2r}}
\le {1.3\over \sqrt{x}} {\rm ~for~}x\ge 289,
\end{equation}
and that, for every fixed $v>1$ and every $x>0$,
\begin{equation}
\label{flauwschatting}
{\log v\over v-1}(1-{v\over x})\le \sum_{r=1}^{[\log x/\log v]}{\log v\over v^r}
\le {\log v\over v-1}.
\end{equation}
\begin{Thm} 
\label{theorem3}
We have\\
{\rm a)} $C_{-}(b_1)=-\log \sqrt{2}$ and $C_{+}(b_1)<0.2663$\\
{\rm b)} $C_{-}(b_3)=-\log \sqrt{3}$ and $C_{+}(b_3)<0.276$
\end{Thm}
{\it Proof}. After some computation for the interval $[1,68]$ we
infer, from Theorem \ref{olivier}, that $C_{+}(b_1)\le B_{b_1}+0.0961$
and similarly $C_{+}(b_3)\le B_{b_3}+0.1205$. For 
the determination of $C_{-}(b_i)$ we use
(\ref{tochnodig}) and
(\ref{flauwschatting}) in addition to Ramar\'e's inequalities,
which yields, for $x\ge 289$, that
$H_1(x)\ge {1\over 2}\log x+B_{b_1}-{1.3/\sqrt{x}}
-{(\log 4)/x}$ and  
$H_3(x)\ge {1\over 2}\log x+B_{b_3}-{1.3/\sqrt{x}}
-(\log 27)/(2x)$.\qed\\

\noindent Let RH$(d)$ be the hypothesis that for every character $\chi$ mod $d$
every non-trivial zero of $L(s,\chi)$ is on the critical line.
\begin{Thm} We have\\
\label{supinf}
{\rm a)} $C_{+}(b_1)=H_1(461)=0.1701069880305239\cdots$, under RH$(4)$.\\
{\rm b)} $C_{+}(b_3)=H_3(3739)=0.1554480047272349\cdots$, under RH$(3)$.
\end{Thm}
{\it Proof}. (cf. \cite[Theorem 6]{moree}).
We recall from \cite{moree} that
for $d\le 432$ and $(a,d)=1$, there exists a constant $c_{d,a}$ such that
for $x\ge 224$ we have, on RH$(d)$, that
$$\left|\sum_{n\le x\atop n\equiv a({\rm mod~}d)}
{\Lambda(n)\over n}-{\log x\over \varphi(d)}-c_{d,a}\right|
\le {11\over 32\pi\sqrt{x}}\{3\log^2 x+8\log x+16\}.$$
Under RH(4) it follows from this that 
$C_{+}(b_1)=\max_{v_i\le 6.15\times 10^8}H_1(v_i)$, where $2=v_1<v_2<\cdots$ are
the consecutive prime powers that can be written as a sum of two squares.
Similarly under RH(3) we deduce that
$C_{+}(b_3)=\max_{w_i\le 1.083\times 10^{10}}H_3(w_i)$, where $3=w_1<w_2<\cdots$ 
are
the consecutive prime powers that can be represented by the form
$X^2+3Y^2$. On computing these maxima (for details
see Section \ref{herman}), the proof is then
completed. \qed\\

\noindent The reason that, even on GRH, it requires a lot of computation to 
determine
$C_{+}(b_1)$ and $C_{+}(b_3)$ is that these values are so close to
$B_{b_1}$, respectively $B_{b_3}$. A similar phenomenon occurs in
\cite{moree} for some of the functions considered there (cf. Theorem 6).\\
\indent Using Theorem \ref{theorem3} and Lemma \ref{lemma1} together
with good enough approximations for $C_{b_1}$ and $C_{b_3}$, one infers
that $\mu_{b_1}(x)\ge \mu_{b_3}(x)$ for $x\ge 27500$. After some
computation we then deduce that Conjecture 3 holds true.\\
\indent Unfortunately
establishing Conjecture 2 requires quite a bit more work. In particular
we need values for $D_{-}$ and $D_{+}$ in Lemma \ref{lemma1} that are
more closely together than those coming from Theorem \ref{theorem3}.
Without improvement of Theorem \ref{olivier}, the upper bounds in 
Theorem \ref{theorem3} cannot be improved.  The lower bounds, however,
are amenable to improvement.\\
\indent Let $\Delta_f(x)$ denote the quantity that is sandwiched between
$D_{-}\mu_f(x)$ and $D_{+}\mu_f(x)$ in
(\ref{volgendeinklemming}).
Using the lower bound for $H_3(x)$ appearing in
the proof of Theorem \ref{theorem3}, we deduce
that $H_3(x)\ge 0$ for $x\ge 25$.
We infer that 
$$\Delta_{b_3}(x)\ge
-{\log \sqrt{3}}\{\mu_{b_3}(x)-\mu_{b_3}({x\over 25})\}.$$
On applying Lemma \ref{delaatstehoopik} with
$D_{-}=-\log \sqrt{3}$ and $D_{+}=0.276$, we deduce that
$\Delta_{b_3}(x)/\mu_{b_3}(x)\ge -0.09586\cdots$ for
$x\ge 10^9$. Taking  $D_{-}=-0.09586\cdots$ as new value
and repeating the procedure, we obtain $D_{-}=-0.06890\cdots$.
Iterating twice more, we see that 
for $x\ge 10^9$ we can take $D_{-}=-0.0672$ in Lemma \ref{lemma1}.\\
\indent For any $x$ satisfying the conditions of Lemma \ref{lemma1}, we can
proceed as above. If the first iteration yields an improved value
of our initial $D_{-}$ (which we take to be $-\log \sqrt{3}$), then it is not 
difficult to see that every further iteration
yields a value of $D_{-}$ not less than the previous one (this
is so since, for given $r\ge 1$, 
$L(x/r,{1\over 2},D_{-},D_{+})/U(x,{1\over 2},D_{-},D_{+})$ is
increasing, considered as a function $D_{-}$). On the
other hand the value cannot be improved beyond zero and hence the
iteration process must converge. 
If the first iteration does not yield an improved value
for $D_{-}$ (which is initially
taken as $-\log \sqrt{3}$), we put ${\tilde w}_i(x)=-\log \sqrt{3}$
for every $i\ge 0$, otherwise we put 
${\tilde w_0}(x)=-\log \sqrt{3}$ and define 
$${\tilde w}_{i+1}(x)=
\left({L({x\over 25}.{1\over 2},{\tilde w}_i(x),0.276)\over U(x,{1\over 
2},{\tilde w}_i(x),0.276)}-1\right){\log 3\over 2}.$$
Empirically it seems
that  
after $n$ iterations we can expect to have approached the limit value
$\lim_{i\rightarrow \infty}{\tilde w}_i(x)$ with $O(n)$ decimal precision.\\
\indent For $b_1$ we proceed
similarly. After some computation using the lower bound
for $H_1(x)$ given in Theorem \ref{theorem3}, we find that
$H_1(x)\ge 0.065$ for $x\ge 97$. Hence
$$H_1(x)\ge
\left(-\log \sqrt{2}-0.065\right)\{\mu_{b_1}(x)-\mu_{b_1}({x\over 
97})\}+0.065\mu_{b_1}(x).$$
If the first iteration does not yield an improved value
for $D_{-}$ (which is initially
taken as $-\log \sqrt{2}$), we put ${\tilde v}_i(x)=-\log \sqrt{2}$
for every $i\ge 0$,
otherwise we put
${\tilde v}_0(x)=-\log \sqrt{2}$ and define 
$${\tilde v}_{i+1}(x)=\left({L({x\over 97}.{1\over 2},{\tilde v}_i(x),0.2663)
\over U(x,{1\over 2},{\tilde v}_i(x),0.2663)}-1\right)(\log 
\sqrt{2}+0.065)+0.065.$$
To sum up, we have established:
\begin{Lem}
\label{spelhetmaaruithoor}
Suppose that $x\ge x_0\ge 2$ and $i\ge 0$. Then {\rm (\ref{volgendeinklemming})} 
holds true
with $f=b_1$, $D_{-}={\tilde v}_i(x_0)$, $D_{+}=0.276$. It also
holds true with $f=b_3$, $D_{-}={\tilde w}_i(x_0)$ and $D_{+}=0.2663$.
\end{Lem}
This lemma, although amenable to
further improvement, is sufficiently sharp for our purposes.

\section{The proof of Theorem \ref{hiergaathetom}}
Before proving Theorem \ref{hiergaathetom},
we will need two more lemmas. From prime number theory we
recall that
$\psi(x;d,a)=\sum_{n\le x,~n\equiv a({\rm mod~}d)}\Lambda(n)$.
\begin{Lem}
\label{psiafschat} We have\\
{\rm a)} $\psi_{b_1}(x)\ge 0.4924x$ for $x\ge 37$.\\
{\rm b)} $\psi_{b_3}(x)\le 0.5176x$ for $x\ge 3793$.
\end{Lem}
{\it Proof}. Let $d\le 13$ and $(a,d)=1$. Then 
$|\psi(x;d,a)-x/\varphi(d)|\le \sqrt{x}$ 
for $224\le x\le 10^{10}$
by \cite[Theorem 1]{ramare}
and $|\psi(x;d,a)-{x/\varphi(d)}|<0.004560{x/\varphi(d)}$ for $x\ge 10^{10}$
by \cite[Theorem 5.2.1]{ramare}.
From these inequalities the lemma 
follows after some computation. \qed \\

\noindent For $y\ge 3$ we define $S_{b_3}(y)$ by $0.5176y$, except for the
intervals $[3,49),[49,181)$, $[181,487),[487,1369),[1699,1933),[2287.2437)$,
respectively $[3733,3793)$, where we define $S_{b_3}(y)$ to be
(respectively) $0.653954y,0.605778y,0.557372y,0.534528y$, $0.526579y,0.521825y$
and $0.51996y$.
\begin{Lem}
\label{shapetracer}
For $y\ge 2$ we have $\psi_{b_3}(y)\le S_{b_3}(y)$.
\end{Lem}
{\it Proof}. The points where $\psi_{b_3}$ and $S_{b_3}$ change value
occur only at prime powers representable by $X^2+3Y^2$, which we
denoted by $3=w_1<w_2<\cdots$.
We now
check that $\psi_{b_3}(w_i)\le S_{b_3}(w_i)$ for every $w_i\le 3793$.
For $w_i\ge 3793$ the result follows by Lemma \ref{psiafschat}. \qed
\\ \\
${~~~~~~~~}${\resizebox{5.0cm}{!}{\includegraphics{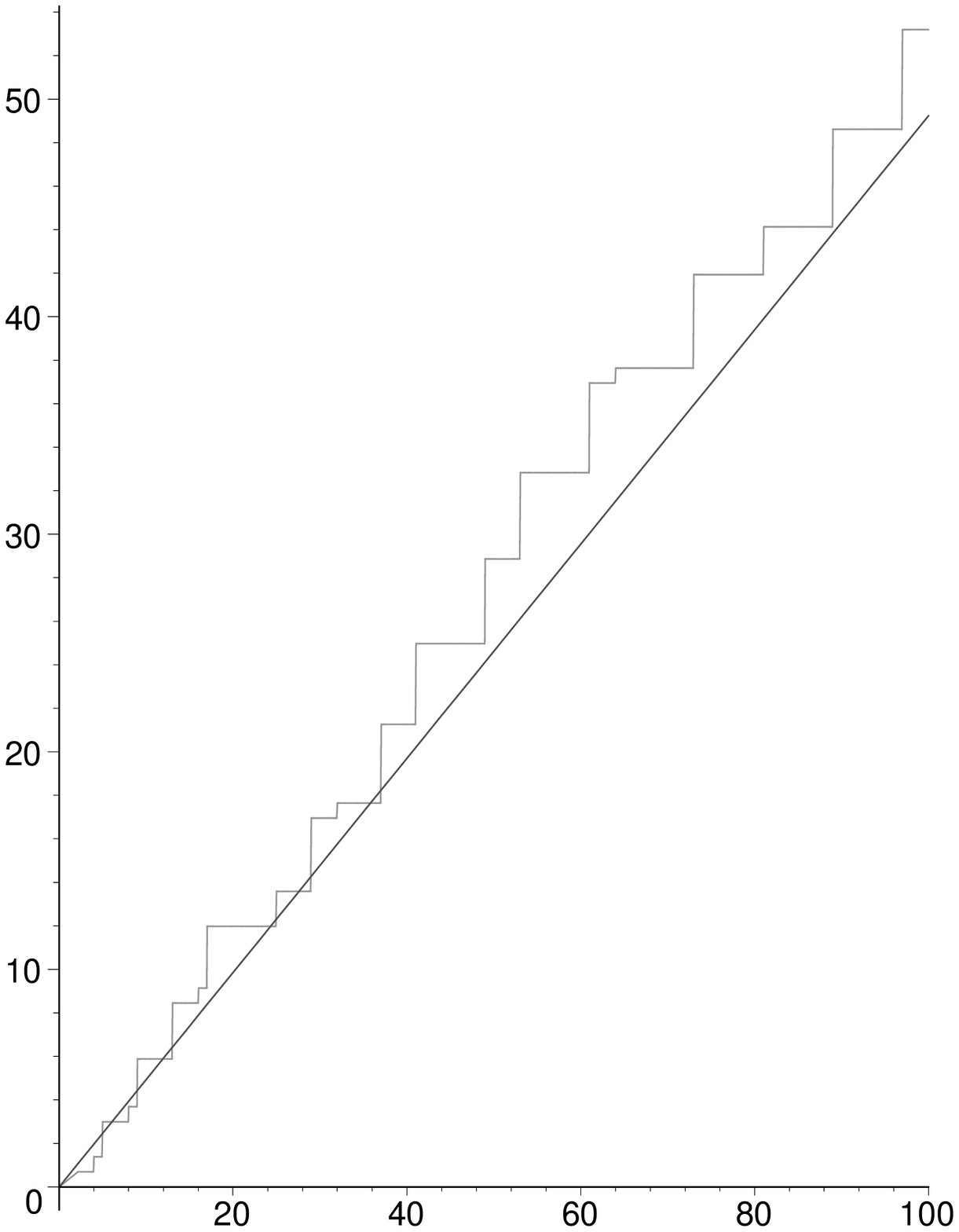}}}
${~~~~~~~~~~~~~}${\resizebox{5.0cm}{!}{\includegraphics{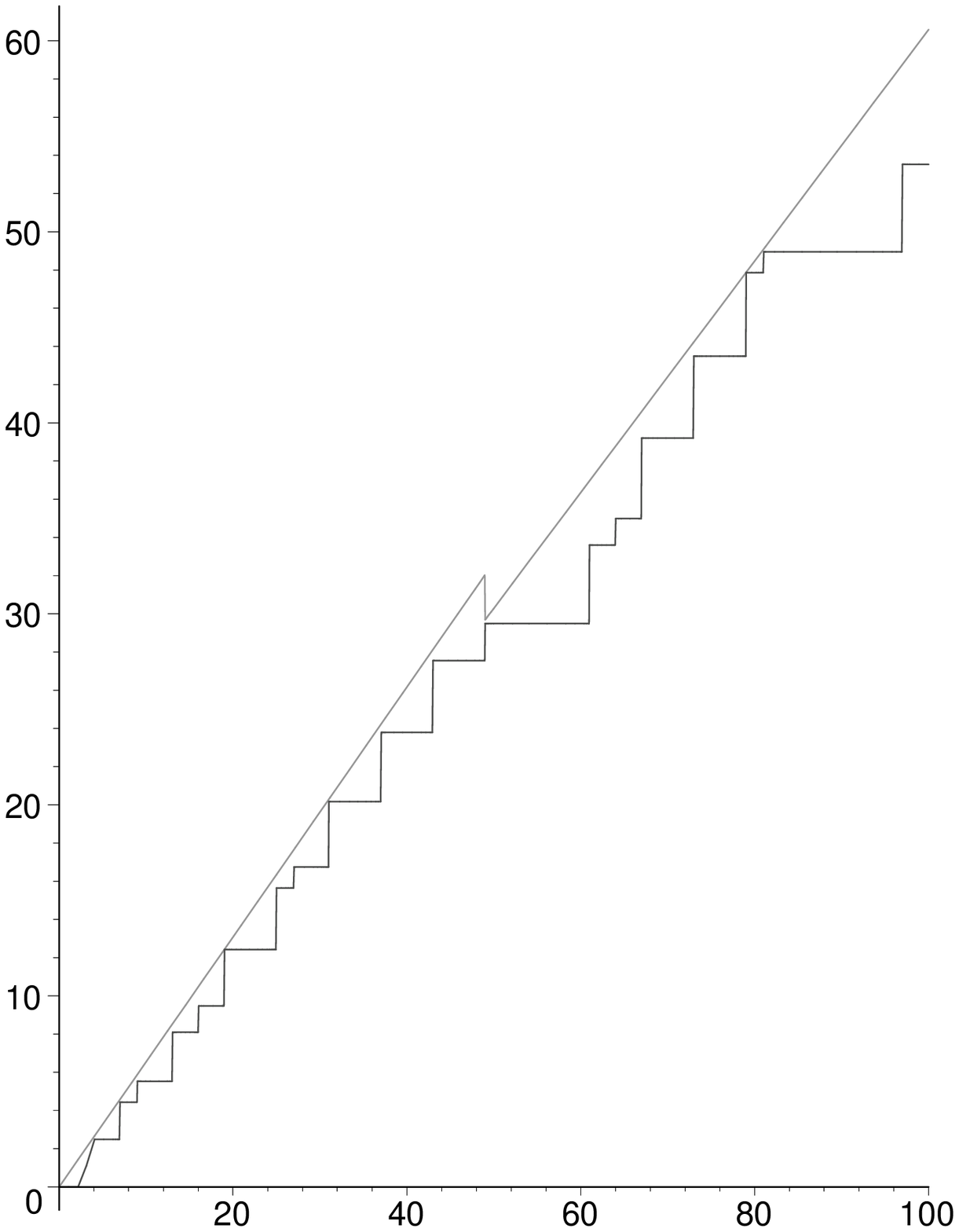}}}\\
${~~~~~~~~~~~~~~~~~~~~~~}${\bf Figure 3}
${~~~~~~~~~~~~~~~~~~~~~~~~~~~~~~~~~~~~}$ {\bf Figure 4}\\
${~~~~~~~~}${\tt Plot of $\psi_{b_1}(x)$ versus $0.4924x$}
${~~~~~~}${\tt Plot of $\psi_{b_3}(x)$ versus $S_{b_3}(x)$}\\
${~~~~~~~~~~~~~~~~~~}$ {\tt for $0\le x\le 100$}
${~~~~~~~~~~~~~~~~~~~~~~~~~~~}${\tt for $0\le x\le 100$}
\\ \\
%
%
%
%
\indent At last we are in the position to prove Theorem \ref{hiergaathetom}.\\

\noindent {\it Proof of Theorem} \ref{hiergaathetom}. As we have 
shown in Section 2, it suffices to establish Conjecture 2. To
this end we have to prove that, for $x\ge 8$,
$$\lambda_{b_1}(x)=\sum_{n\le {x\over 2}}b_1(n)\psi_{b_1}({x\over n})\ge 
\sum_{n\le {x\over 3}}b_3(n)\psi_{b_3}({x\over n})=\lambda_{b_3}(x).$$
Let us denote the 6 intervals in the definition of $S_{b_3}(y)$ by
$[r_i,s_i)$ for $i=1,\cdots,6$ and put $\alpha_i=S_{b_3}(r_i)/r_i-0.5176$
(note that $\alpha_i>0$).
From Lemma \ref{shapetracer} we infer that
$$\lambda_{b_3}(x)\le \sum_{n\le {x\over 3}}b_3(n)S_{b_3}({x\over n})
=0.5176\mu_{b_3}({x\over 3})+\sum_{i=1}^6 \alpha_i\{\mu_{b_3}({x\over r_i})-
\mu_{b_3}({x\over s_i})\}.$$
Put $x_0=1.5\times 10^{11}$. 
Using a computer (see Section \ref{herman}) Conjecture 2 can be established for 
$x<x_0$.
Hence now assume that $x\ge x_0$. 
For 
notational convenience
we shorten $U(x/r,{1\over 2}, {\tilde w}_8(x_0/r), 0.276)$ to
$U_3(x/r)$, $L(x/r,{1\over 2}, {\tilde w}_8(x_0/r), 0.276)$ to
$L_3(x/r)$ and $L(x/r,{1\over 2}, {\tilde v}_8(x_0/r), 0.2663)$
to $L_1(x/r)$, where $r$ is some fixed number.
On applying Lemma \ref{spelhetmaaruithoor}, we deduce that
$${\lambda_{b_3}(x)\over 2C_{b_3}}\le 0.5176U_3({x\over 3})+\sum_{i=1}^6 
\alpha_i\{U_3({x\over r_i})-
L_3({x\over s_i})\}.$$
By Lemma \ref{dalertje} each of
the six terms in the above sum is non-increasing for
$x\ge x_0$ and thus the sum is bounded above by its value
in $x_0$, which on its turn is less 
than $0.0224U_{3}(x_0/3)$. One easily checks that
$U_{3}(x/3)\ge U_{3}(x_0/3)$ for $x\ge x_0$ (on noting that
$U(y,\tau,D_{-},D_{+})$, considered as a function of
$\tau$, is increasing for $y>\exp(D_{+}+(D_{+}-D_{-})/\tau)$).
We thus obtain that $\lambda_{b_3}(x)
\le 1.08C_{b_3}U_{3}(x/3)$.
Using Lemma \ref{spelhetmaaruithoor} 
and the lower bound for $\psi_{b_1}$ given
in Lemma \ref{psiafschat}, we infer
that
$\lambda_{b_1}(x)\ge \sum_{n\le x/37}b_1(n)\psi_{b_1}(x/n)
\ge 0.4924\mu_{b_1}(x/37)\ge 0.9848C_{b_1}L_1(x/37)$. 
A computation shows that $0.9848C_{b_1}L_{1}(x_0/37)>
1.08C_{b_3}U_{3}(x_0/3)$. By
Lemma \ref{delaatstehoopik} we then have 
$0.9848C_{b_1}L_{b_1}(x/37)>
1.08C_{b_3}U_{3}(x/3)$. 
 for every $x\ge x_0$.
We thus obtain that for every $x\ge x_0$,
$$\lambda_{b_1}(x)\ge 0.9848C_{b_1}L_{1}({x\over 37})
\ge 1.08C_{b_3}U_{3}({x\over 3})\ge \lambda_{b_3}(x),$$
completing the proof. \qed

\section{Computations of results used in Theorems \ref{hiergaathetom} and 
\ref{supinf}}
\label{herman}
In the proof of Theorem \ref{hiergaathetom} we have used that Conjecture 2 is 
true
for $x\le x_0$ with $x_0=1.5\times10^{11}$. We have established that result
as follows.
\par
Checking Conjecture 2 requires the computation and comparison of the sums 
$$\lambda_{b_i}(x)=\sum_{n\le x}b_i(n)\log n,~~i=1,3,$$
and, consequently, the computation of 
the characteristic functions $b_1(n)$ and $b_3(n)$ for all positive integers
$n\le x_0$. 
Because of the size of $x_0$, the range of $x$-values for which Conjecture 2 
had to be checked, was split up
in sub-intervals of length $10^6$, large enough for efficiency, and small
enough to avoid so-called cache misses during the computations. 
\par
We first describe the case $\lambda_{b_1}(n)$.
For a given interval, say, $[A,B]$, an integer array $b(j),j=1,2,\dots,B-A+1$
of length $B-A+1$ is initialized to $0$. Here, $b(j)$ corresponds to 
$b_1(j+A-1)$. 
Next, all the possible sums of squares $x^2+y^2$ of integers $0\le x\le y$,
with $A\le x^2+y^2\le B$, hence $x\in\left[0,\sqrt{B/2}\right]$, 
$y\in\left[\sqrt{A/2},\sqrt{B}\right]$,
are computed as follows. First, the sequence of all the squares 
$y^2\in[A/2,B]$ is precomputed and stored. 
Next, for each $x=0,1,\dots,\left\lfloor\sqrt{B/2}\right\rfloor$, 
the sums $x^2+y^2$ are computed 
$$\mbox{for all~~} y^2\in\left[\max(A-x^2,A/2),B-x^2\right]\subset[A/2,B].$$ 
For all the sums $x^2+y^2=:n$ obtained in this way, $b(n-A+1)$ is set to $1$.
\par
The case $\lambda_{b_3}(n)$ is treated similarly: the same initialization of 
array $b$ is carried out. Next, all the possible sums $x^2+3y^2$ of integers $x, 
y$, with $A\le x^2+3y^2\le B$, hence
$x\in\left[0,\sqrt{B}\right]$, $y\in\left[0,\sqrt{B/3}\right]$
are computed as follows. First, the sequence of all the triples of squares
$3y^2\in[0,B]$ is precomputed and stored. 
Next, for each $x=0,1,\dots,\left\lfloor\sqrt{B}\right\rfloor$,
the sums $x^2+3y^2$ are computed 
$$\mbox{for all~~} 3y^2\in\left[\max(0,A-x^2),B-x^2\right]\subset[0,B]$$ 
and for all the sums $x^2+3y^2=:n$ obtained in this way, 
$b(n-A+1)$ is set to $1$. This corresponds to $b_3(n)$. 
\par
We have implemented these algorithms for $b_1(n)$ and $b_3(n)$
in Fortran and used them to compute $\lambda_{b_i}(x)$ for $i=1, 3$,
and to verify Conjecture 2 for $x=8,9,\dots,1.5\times10^{11}$
on one 250 MHZ processor of CWI's SGI Origin 2000 computing system. Computing
time was 7.6 CPU hours. We also used our program to check the values 
of $B_1(x)=\sum_{n\le x}b_1(n)$, given for $x=10^i, i=1,\dots,12$,
by Shiu in Table 1 of \cite{shiu}
(where $B_1(x)$ is called $W(x)$). 
Computing time to extend our results from $1.5\times10^{11}$ to $10^{12}$ was
77 CPU hours. 
We found agreement with Shiu 
for $i=1,\dots,10$, but differences for $i=11$ and $i=12$:
$B_1(10^{11})=15~570~512~744$ and $B_1(10^{12})=148~736~628~858$,
whereas Shiu gave $W(10^{11})=15~570~523~346$ and $W(10^{12})=148~736~629~005$.
Shiu used a different, more efficient method than ours, but he has confirmed
our value of $B_1(10^{11})$ after checking and correcting his program 
\cite{shiuprivate}.
\par
We have spot-checked our program for computing $b_1(n)$ and $b_3(n)$ on various 
intervals of length $10^6$
with the help of Lemma \ref{classical}. This requires the decomposition in 
primes of  
each $n$ for which we wish to compute $b_1(n)$, which is 
extremely expensive, compared with composing all integers in a given long 
interval
$[A,B]$ as a sum of integer squares. However, we found agreement for all
the checks we did, in particular for those in the neighboorhood of $x=10^{12}$.
In Table 2, we list, for $i=1,3$, the values we found of 
$\lambda_{b_i}(x)$ and $B_i(x)$ for $x=j\times10^{11}$, $j=1,1.5,2,\cdots,10$.
\begin{table}
\begin{center}
\begin{tabular}{crrrr}\hline
$x/10^{11}$ & $\lambda_{b_1}(x)$ & $\lambda_{b_3}(x)$ & $B_1(x)$ & $B_3(x)$\\
\hline
1 & 378458908590.818 & 316358774044.179 & 15570512744 & 13015595425 \\
1.5 & 572353849423.260 & 478438468735.511 & 23160971166 & 19360573686 \\
2 &   767521856517.400 & 641582406621.494 & 30700929088 & 25663340448 \\
3 & 1160486988190.213 & 970068358550.987 & 45678037444 & 38182949191 \\
4 & 1555965223692.576 & 1300655152892.098 & 60558145064 & 50621477125 \\
5 & 1953301629004.525 & 1632795521743.015 & 75367348255 & 63000746043 \\
6 & 2352112868630.901 & 1966168966371.294 & 90120785046 & 75333407591 \\
7 & 2752146230205.959 & 2300563843364.554 & 104828319151 & 87627692348 \\
8 & 3153223047545.408 & 2635831188875.970 & 119496904413 & 99889427349 \\
9 & 3555209733889.339 & 2971859287714.156 & 134131682979 & 112122909167 \\
10 & 3958003171956.632 & 3308561817015.470 & 148736628858 & 124331455166 \\
\hline
\end{tabular}
\caption{$\lambda_{b_1}(x), B_1(x)$ versus $\lambda_{b_3}(x), B_3(x)$}
\end{center}
\end{table}
\par
In the proof of Theorems 5a and 5b, we have used that 
\begin{equation}
\max_{v_i\le6.15\times10^8}H_1(v_i)=H_1(461)=0.170106...,,{\rm ~respectively},
\label{Thm5a}
\end{equation}
\begin{equation}
\max_{w_i\le1.083\times10^{10}}H_3(w_i)=H_3(3739)=0.155448...,
\label{Thm5b}
\end{equation}
We have established these results as follows.
\par 
Let $x=6.15\times10^8$. We first generated the primes $\le \sqrt{x}$ with
the sieve of Eratosthenes, 
and stored the following pairs $(n,\Lambda_{b_1}(n))$: 
$$(2^k,\log 2), k=1,2,\dots,\lfloor\log_2x\rfloor,$$
$$(p^{2k},2\log p), k=1,2,\dots,\left\lfloor\frac12\log_px\right\rfloor, 
\mbox{~for~the~primes~}p\equiv3\bmod{4}\le\sqrt{x},$$
$$(p^k,\log p), k=1,2,\dots,\lfloor\log_px\rfloor, 
\mbox{~for~the~primes~}p\equiv1\bmod{4}\le\sqrt{x},$$
into an array, sorted increasingly according to the first element of the pairs. 
The set of numbers $n$ in these pairs in fact contains 
as a subset all the prime powers 
$v_1, v_2, \dots \le \sqrt{x}$ which can be written as a sum of two squares.
For these $(n,\Lambda_{b_1}(n))$-pairs, we computed $H_1(n)$ and we verified 
that 
$$\max_{v_i\le\left\lfloor\sqrt{6.15\times10^8}\right\rfloor}H_1(v_i)=H_1(461)=0
.170106...$$
The remaining interval 
$\left[\left\lfloor\sqrt{6.15\times10^8}\right\rfloor+1,x\right]$ was split up 
in pieces of length $10^7$, and for each of these intervals, $[A,B]$, say, 
the primes $p\equiv1\bmod{4}$ were generated 
with the sieve of Eratosthenes, together with $\log p$. These 
pairs $(p, \log p)$ were mixed with the 
$(n,\Lambda_{b_1}(n))$-pairs generated above for which $n\in[A,B]$
and then it was verified that
$\max_{v_i\in[A,B]}H_1(v_i)<H_1(461)$.
This proved (\ref{Thm5a}). Computing time was 81 CPU seconds. 
Relation (\ref{Thm5b}) was proved in a similar way at the expense of 1340 CPU 
seconds.

\section{An alternative approach}
In the previous sections we have made essential use of asymptotic
information regarding the distribution of primes. Some of the results
we used depend eventually on RH$(3)$ and RH$(4)$ to be true up to
some finite height. It might come as a surprise then that it
is possible to show that $B_1(x)\ge B_3(x)$ for $x\ge 10^{9111}$, without
invoking any result from
computational prime number theory (one only needs the ability to compute some
successive primes...).\\
\indent Our method of establishing this is
inspired on Selberg's 
\cite[pp. 183-185]{selberg} method of obtaining an asymptotic evaluation
for $N(x;4,1)$, where
$N(x;d,a)$ denotes the number of integers $n\le x$ that have no
prime factor $p$ with $p\not\equiv a({\rm mod~}d)$.
 Unfortunately Selberg's method
does not seem to generalise well; for example we
have no idea how to generalise it so as to show that $N(x;4,3)\ge N(x;4,1)$
for $x\ge x_0$, with $x_0$ some effectively computable constant.
Vide \cite{saradha}
for generalisations of Selberg's method.
\begin{Lem}
\label{explicit}
\label{selberglower} a) For $x\ge 2$ we have
$$|B_1(x)-C_{b_1}{x\over \sqrt{\log x}}|\le 9.62{x\over \log x},$$
b) For $x\ge 2$ we have
$$|B_3(x)-C_{b_3}{x\over \sqrt{\log x}}|\le 8.53{x\over \log x},$$
\end{Lem}
\begin{cor}
For $x\ge 10^{9111}$ we have $B_1(x)\ge B_3(x)$.
\end{cor}
\indent In the proof of Lemma \ref{explicit} we will make use of the following 
result.
\begin{Lem}
\label{eenovern}
Let $c_2=2e^{\gamma}$ and $c_3=\sqrt{3}e^{\gamma}$.
For $z\ge 1$ we put
$$f(z)=z\sum_{n\le z,~2\nmid n}
{1\over n}-{z\over 2}\log(c_2z){\rm ~and~}
g(z)=z\sum_{n\le z,~3\nmid n}
{1\over n}-{2\over 3}z\log(c_3z).
$$
We have
$$\sup_{z\ge 1}|f(z)|=-f(3^{-})={3\over 2}\{\log 6+\gamma\}-3
=0.55346270119438\cdots$$
and
$$\sup_{z\ge 1}|g(z)|=-g(4^{-})={8\over 3}\log(4\sqrt{3}e^{\gamma})-6=
0.70084312094794\cdots$$
\end{Lem}
{\it Proof}. We only prove the statement concerning $f(z)$, the 
statement regarding $g(z)$ can be proved in a similar way.\\ 
\indent Since $z\log(c_2z)$ is monotonically increasing, we obtain,
cf. \cite[Lemma 4]{moree}, that
$\sup_{z\ge 1}|f(z)|=\sup\{|f(1)|,|f(3^{-})|,|f(3)|,|f(5^{-})|,|f(5)|,
\cdots \}$.  
Using the Euler-MacLaurin summation formula, 
cf. \cite[p. 6]{tenenbaum}, one finds that for integers
$n\ge 1$,
\begin{equation}
\label{maclaurin}
\sum_{m\le n}{1\over m}=\log n+\gamma+{1\over 2n}-{1\over 12n^2}
+{\theta_1(n)\over 60n^4},
\end{equation}
where $\theta_1(n)\in [0,1]$.
Clearly
\begin{equation}
\label{maclaurin2}
\sum_{m\le n,~2\nmid m}{1\over m}
=\sum_{m\le n}{1\over m}-{1\over 2}\sum_{m\le n/2}{1\over m},
\end{equation}
Let $n\ge 3$ be an odd integer. Notice that
$f(n^{-})=f(n)-1$. Using (\ref{maclaurin}) and
(\ref{maclaurin2}) it is not difficult to deduce that
$$-{1\over 2}\le f(n)\le {n\over 2}\log({n\over n-1})+{n\over 6(n-1)^2}+{1\over 
60n^3}$$
and
$$-{1\over 2}\le 
-f(n^{-})\le {1\over 2}+{1\over 2(n-1)}+{1\over 12n^2}+{2\over 15(n-1)^4}.$$
Using that the latter two right hand
sides are monotonically decreasing in
$n$ and noting that $\sup_{z\ge 1}|f(z)|\ge |f(3^{-})|$, we
see
 that
$\sup_{z\ge 1}|f(z)|=\sup_{1\le z\le 11}|f(z)|=|f(3^{-})|$.
\qed\\

\noindent Let
$\{z\}=z-[z]$ denote
the fractional part of $z$. Using (\ref{maclaurin}) it can be shown that the functions
$f$ and $g$ are almost periodic in the sense that they converge uniformly
to the periodic functions ${1\over 2}-\{ {z-1\over 2}\}$, respectively
$1-\{ {z-1\over 3}\}-\{{z+1\over 3}\},$ 
 (cf. Figure 5).
\\ \\
${~~~~~~~~}${\resizebox{5.0cm}{!}{\includegraphics{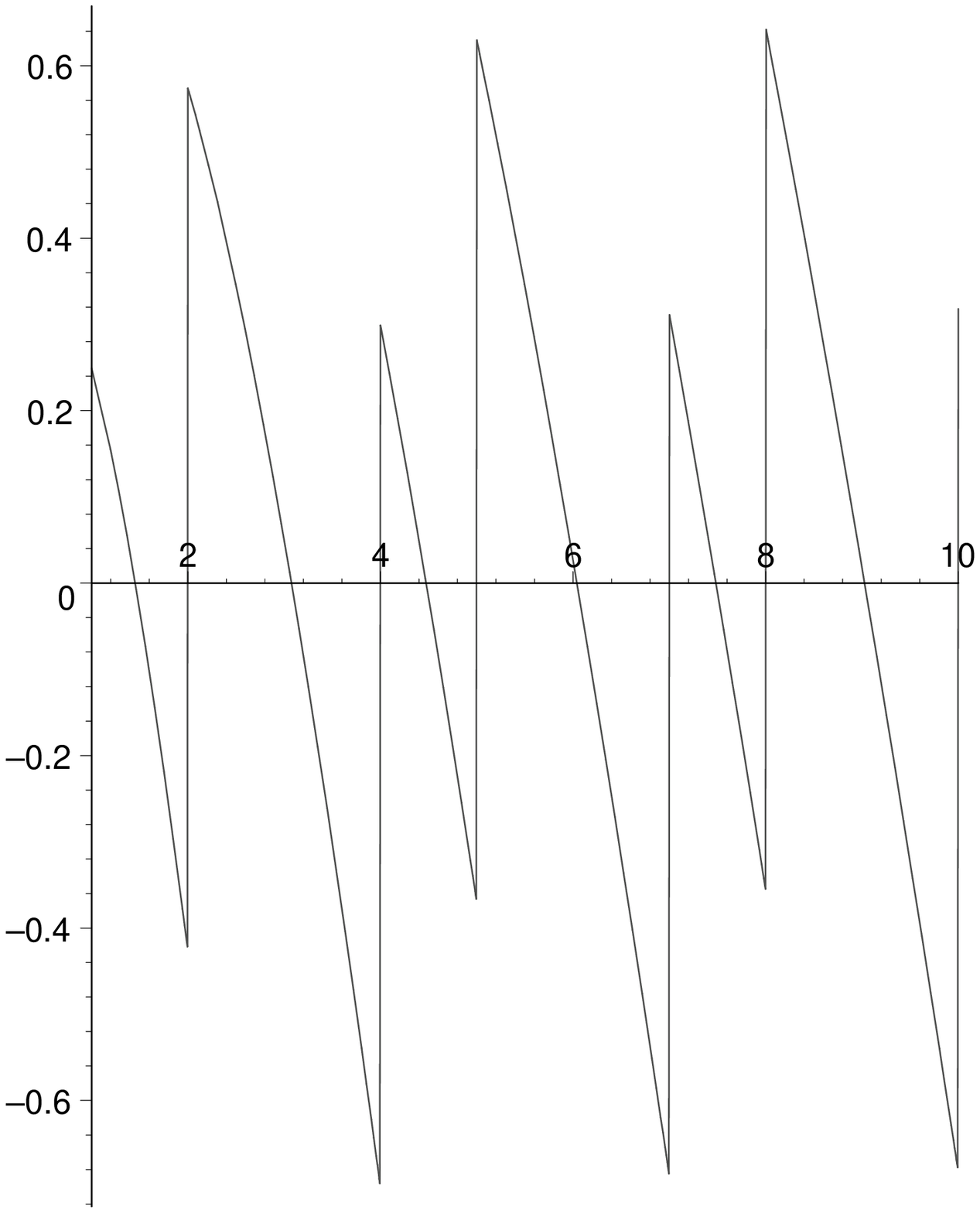}}}
${~~~~~~~~~~~~~}${\resizebox{5.0cm}{!}{\includegraphics{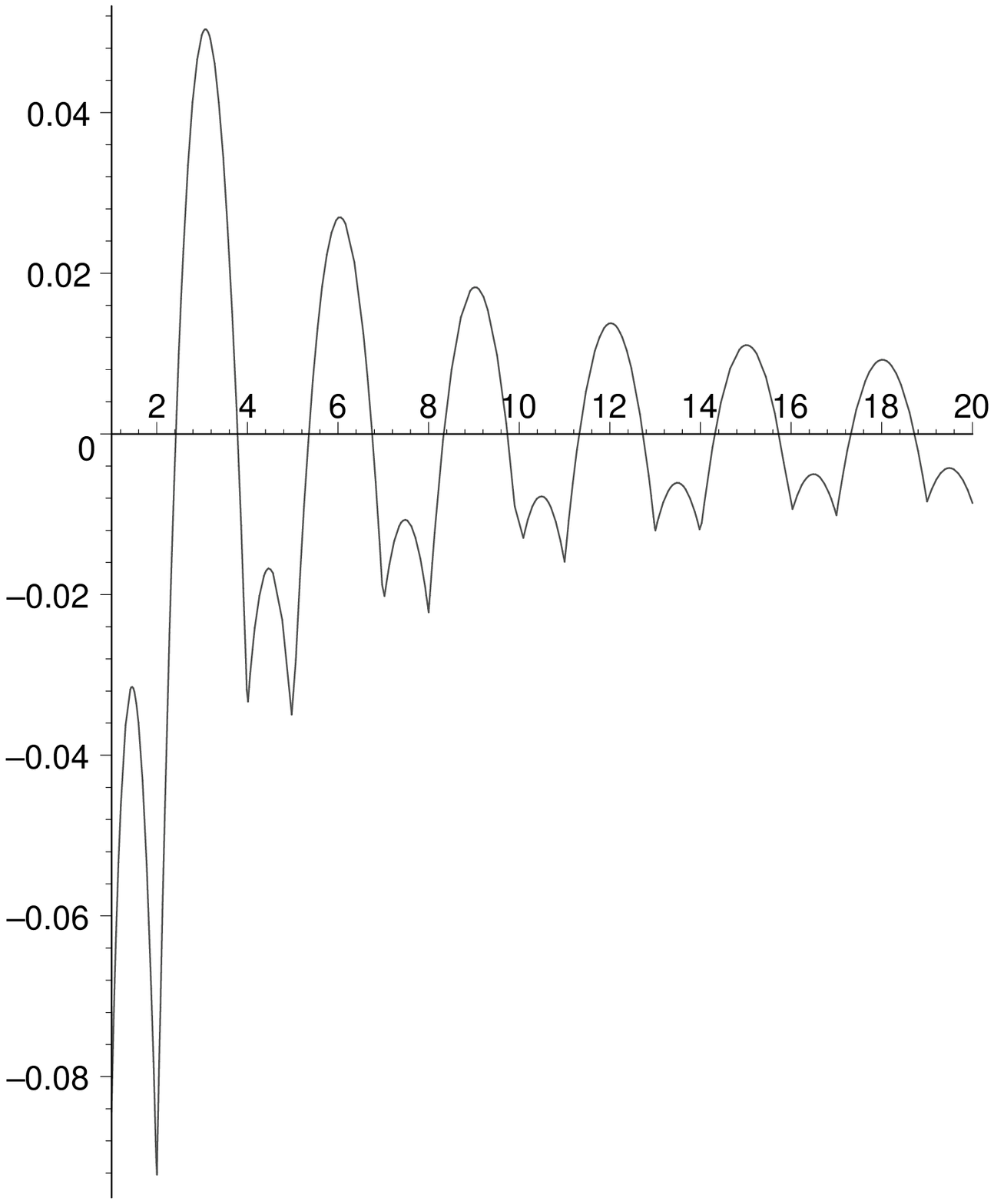}}}\\
${~~~~~~~~~~~~~~~~~~~~~~}${\bf Figure 5}
${~~~~~~~~~~~~~~~~~~~~~~~~~~~~~~~~~~~~}$ {\bf Figure 6}\\
${~~~~~~~}${\tt Plot of $g(x)$ for $1\le x\le 10$}
${~~~~~}${\tt Plot of $g(x)$ minus its limit function}\\
${~~~~~~~~~~~~~~~~~~~~~~~~~~~~~~~~~~~~~~~~~~~~~~~~~~~~~~~~~~~~~~~~~~}$ {\tt for 
$1\le x\le 20$}
%
\\
\\
\section{Proof of Lemma \ref{selberglower}}
Let $P_2, P_3$ denote the 
set of primes $p$ that satisfy $p\equiv 2({\rm mod~}3)$,
respectively $p\equiv 3({\rm mod~}4)$. Let $(P_2),(P_3)$ denote the set
of natural numbers that have no 
prime divisor $p$ with $p\not\equiv 2({\rm mod~}3)$, respectively 
$p\not\equiv 3({\rm mod~}4)$.
Let $\psi_3(x),\psi_4(x)$ denote the number of integers $1\le n\le x$ 
that have no prime divisor $p$ with $p\not\equiv 2({\rm mod~}3)$, respectively 
$p\not\equiv 3({\rm mod~}4)$.\\ 

\noindent {\it Proof of part a}. Put $c_2=2e^{\gamma}$.
We consider the expression
$$
\sum_{j=0}^{\infty}
\sum_{m\in (P_3)}\sum_{1\le n\le x/(2^jm^2)\atop n\equiv 1({\rm 
mod~}4)}\sum_{d|n\atop
d\in (P_3)}\mu(d)\log{c_2x\over 3^jm^2d}=$$
\begin{equation}
\label{basis}
\sum_{j=0}^{\infty}
\sum_{m\in (P_3)}
\sum_{d\in (P_3)\atop d\le x/(2^jm^2)}
\mu(d)\log{c_2x\over 2^jm^2d}\sum_{d|n,~1\le n\le x/(2^jm^2)\atop n\equiv 1({\rm 
mod~}4)}1.
\end{equation}
By approximating both sides of this equation in terms of 
the function $B_1$, we
will arrive at an approximate functional equation,
(\ref{functional}), for $B_1$ which
on solving will yield an explicit lower bound for $B_1$.\\
\indent
On recalling that $\sum_{d|n}\mu(d)\log d=-\Lambda(n)$, one sees that,
when $n\equiv 1({\rm mod~}4)$ and $d\le x/(2^jm^2)$, we have
$$
\sum_{d|n\atop d\in (P_3)}\mu(d)\log{c_2x\over 2^jm^2d}=
\cases{\log {c_2x\over 2^jm^2} &if $n\in (P_1)$;\cr
\log p &if $n$ is divisible by exactly one $p$ from $P_3$;\cr
0 &if $n$ has $\ge 2$ distinct prime factors from $P_3$.\cr}
$$
On noting that
$$B_1(x)=\sum_{j=0}^{\infty}\sum_{m\in (P_3)}\psi_4({x\over 2^jm^2}),$$
we see that 
the left hand side of (\ref{basis}) equals
\begin{equation}
\label{lefthandside}
B_1(x)\log c_2 x+\sum_{p\in P_3\atop r\ge 1}\log p~B_1({x\over p^{2r}})
-\sum_{j=0}^{\infty}\sum_{m\in (P_3)}\log(2^jm^2)\psi_4({x\over 2^jm^2}),
\end{equation}
which we write as
\begin{equation}
\label{rewritten}
B_1(x)\log c_2 x+I_1(x)-I_2(x).
\end{equation}
We write $z=x/2^j$ and consider the expression formed by the three
inner sums in (\ref{basis}), that is
\begin{equation}
\label{triplesum}
\sum_{m\in (P_3)}
\sum_{d\in (P_3)\atop d\le z/m^2}
\mu(d)\log{c_2 z\over m^2d}\sum_{d|n,~1\le n\le z/m^2\atop n\equiv 1({\rm 
mod~}4)}1.
\end{equation}
Given an integer $k$, let $\sigma_0(k)$ denote the product of the distinct
primes that occur to an odd power in the prime factorisation of $k$.
We put $\sigma_0(k)=1$ if if there is no prime that occurs to an odd
power in $k$.
Note that 
\begin{equation}
\label{simplearith} 
\sum_{m^2d=k}\mu(d)=\mu(\sigma_0(k)),
\end{equation}
 where the sum is over
all integers $m$ and $d$ such that $m^2d=k$. On writing $m^2d=k$ in 
(\ref{triplesum}), and invoking (\ref{simplearith}), we deduce that the
triple sum in (\ref{triplesum}) equals
$$\sum_{k\in (P_3)\atop k\le z}\mu(\sigma_0(k))\log {c_2 z\over k}
\sum_{k_1\le z/k\atop k_1\equiv k ({\rm mod~}4)}1.$$ 
The right hand side of (\ref{basis}) is thus seen to equal
$$\sum_{j=0}^{\infty}\sum_{k\in (P_3)\atop k\le x/2^j}
\mu(\sigma_0(k))\log {c_2 x\over 2^jk}
\sum_{k_1\le x/(2^jk)\atop k_1\equiv k ({\rm mod~}4)}1,$$
which simplifies to
$$\sum_{d\in (P'_3)\atop d\le x}\mu(\sigma_0(d'))\log{c_2 x\over d}
\sum_{k_1\le x/d\atop k_1\equiv d'({\rm mod~}4)}1,$$
where $d'$ denotes the largest odd divisor of $d$ and
$P_3'=P_3\cup \{2\}$ and $(P_3')$ is defined as $(P_3)$ but where
now no prime divisor $p$ with $p\equiv 1({\rm mod~}4)$ is allowed.
The right hand side of (\ref{basis}) is thus seen to equal
\begin{equation}
\label{geenidee}
\sum_{d\in (P_3')\atop d\le x}
\left[{x\over 4d}+{2+(-1)^{d'-1\over 2}\over 4}\right]\mu(\sigma_0(d'))
\log{c_2 x\over d}
={x\over 4}\sum_{d\in (P_3')\atop d\le x}{\mu(\sigma_0(d'))\over d}
\log {c_2 x\over d}+I_3(x),
\end{equation}
where
\begin{equation}
\label{itwee}
|I_3(x)|\le {3\over 4}\sum_{d\in (P_3')\atop d\le x}\log{c_2 x\over d}.
\end{equation}
For $1\le d\le x$
we have, recalling the definition of $f$ (made
in Lemma \ref{eenovern}),
$${1\over 2}\log{c_2 x\over d}=\sum_{n\le x/d,~2\nmid n}
{1\over n}-f({x\over d}){d\over x}.$$
Combining the latter equation with
the sum in the right hand side of (\ref{geenidee}) yields
$${x\over 4}\sum_{d\in (P_3')\atop d\le x}
{\mu(\sigma_0(d'))\over d}\log{c_2 x\over d}=
{x\over 2}\sum_{dn\le x,~d\in (P_3')\atop 2\nmid n}{\mu(\sigma_0(d'))\over dn}+
I_4(x),$$
where
\begin{equation}
\label{i4}
|I_4(x)|\le {1\over 2}\sum_{d\in (P_3')\atop d\le x}|f({x\over d})|.
\end{equation}
Note that
$${x\over 2}\sum_{dn\le x,~d\in (P_3')\atop 2\nmid n}{\mu(\sigma_0(d'))\over 
dn}=
{x\over 2}\sum_{k\le x}{1\over k}\sum_{dn=k,~d\in (P_3')\atop
2\nmid n}\mu(\sigma_0(d')).$$
Denote the latter inner sum by $h(k)$. We claim that $h(k)=b_1(k)$.
First let us consider the case where $k$ is odd. Then
$$h(k)=\sum_{dn=k,~d\in (P_3)}\mu(\sigma_0(d))
=\sum_{dn=k}\mu(\sigma_0(d))2^{-\omega(d)}\prod_{p|d}(1-(-1)^{p-1\over 2}),$$
where $\omega(d)$ denotes the number of distinct primes dividing $d$.
We see that for odd $k$, $h$ is the Dirichlet convolution of two
multiplicative functions and is thus itself a multiplicative function. For
arbitrary $k$ we note that $h(k)=h(k')$, where $k'$ is the largest odd 
divisor of $k$. Thus $h$ is a multiplicative function. An easy 
computation shows that for every prime power $q$ we have $h(q)=b_1(q)$.
Since both $h$ and $b$ are multiplicative, this completes the proof of
the claim. We thus infer that
$$
{x\over 2}\sum_{dn\le x,~d\in (P_3')\atop 2\nmid n}{\mu(\sigma_0(d'))\over dn}
={x\over 2}\sum_{m\le x}
{b_1(m)\over m}={x\over 2}\int_1^x{dB_1(t)\over t}={B_1(x)\over 2}+{x\over 
2}\int_1^x{B_1(t)\over t^2}dt.$$
Thus the right hand side of (\ref{basis}) equals
$${B_1(x)\over 2}+{x\over 2}\int_1^x{B_1(t)\over t^2}dt+
I_3(x)+I_4(x).$$
Equating it with the expression
(\ref{rewritten}) for the left hand side of (\ref{basis}),
we get
\begin{equation}
\label{functional}
B_1(x)\log c_2 x-{x\over 2}\int_1^x{B_1(t)\over t^2}dt=
-I_1(x)+I_2(x)+I_3(x)+I_4(x)+{B_1(x)\over 2}.
\end{equation}
\indent Next we will consider effective estimates for
$I_j(x)$ for $1\le j\le 4$. Using the trivial estimate $B_1(x)\le x$, we obtain 
that
$$0\le I_1(x)\le x\sum_{p\in P_3}{\log p\over p^2-1}
<.23x.$$
On noting that
$$\sum_{j=0}^{\infty}\sum_{m\in (P_3)}{\log(2^jm^2)\over 2^jm^2}=
2\sum_{m\in (P_3)}{\log(2m^2)\over m^2}<2.7$$
and $\psi_4(x)\le x$, we
deduce that $0\le I_2(x)< 2.7x$.
Using (\ref{itwee}) we deduce that
$$|I_3(x)|\le {3\over 4}\sum_{d\in (P_3'),~d\le x}
\log{cx\over d}
\le {3\over 4}\sum_{1\le d\le x}\int_d^{cx}{dt\over t}\le {3\over 4}cx
<2.68x.$$
For $I_4(x)$ we have, by (\ref{i4}) and Lemma \ref{eenovern},
$|I_4(x)|\le 0.277x$.\\
\indent
Put $A(x)=\int_1^x{B_1(t)dt/t^2}$.
An easy calculation (divide by $x^2\log^{3/2}x$ and
integrate) now shows that if
\begin{equation}
\label{fffunctionaal}
-\alpha_{-}x\le x^2\log xA'(x)-{x\over 2}A(x)\le \alpha_{+}x,~{\rm for~}
x\ge x_0,
\end{equation}
then there exists a constant $c_0$ such that
$$c_0\sqrt{\log x}-2\alpha_{+}\le A(x)\le c_0\sqrt{\log x}+2\alpha_{-},~{\rm 
for~}x\ge x_0.$$
On inserting the latter estimate in (\ref{fffunctionaal}) and invoking
(\ref{1908}), it then follows that
\begin{equation}
\label{zeerexpliciet}
\Big|B_1(x)-C_{b_1}{x\over \sqrt{\log x}}\Big|\le (\alpha_{-}+\alpha_{+}){x\over 
\log x},~{\rm for~}
x\ge x_0.
\end{equation}
From our estimates for $I_j(x)$ with $j=1,\cdots,4$,
we see that we can take $\alpha_{-}=3.96$, $\alpha_{+}=5.66$ 
and $x_0=2$. \qed\\

\noindent {\it Proof of part b}. Making the obvious modifications in the 
proof of part a, we deduce that
\begin{equation}
\label{functional2}
B_3(x)\log c_3 x-{x\over 2}\int_1^x{B_3(t)\over t^2}dt=
-J_1(x)+J_2(x)+J_3(x)+J_4(x)+{B_3(x)\over 2},
\end{equation}
where
$J_1(x)=\sum_{p\in P_2,~r\ge 1}\log p~B_3({x\over p^{2r}})$, $J_2(x)=
\sum_{j=0}^{\infty}\sum_{m\in (P_2)}\log(3^jm^2)\psi_3({x\over 3^jm^2}),$
$$|J_3(x)|\le {2\over 3}\sum_{d\in (P_2'),~d\le x}\log{c_3 x\over d}
{\rm ~and~}
|J_4(x)|\le {1\over 2}\sum_{d\in (P_2'),~d\le x}|g({x\over d})|,$$
with
$P_2'=P_2\cup \{3\}$ and $(P_2')$ defined as $(P_2)$, but 
where
now no prime divisor $p$ with $p\equiv 1({\rm mod~}3)$ is allowed.
Reasoning as before we find
$$0\le J_1(x)<x\sum_{p\in P_2}{\log p\over p^2-1}<0.36x{\rm ~and~}
0\le J_2(x)< {3\over 4}x\sum_{m\in (P_2)}{\log(3m^4)\over m^2}<2.7x.$$
Furthermore we find that $|J_3(x)|\le 2xc_3/3<2.06x$ and $|J_4(x)|\le 0.36x$.
From these estimates and (\ref{functional2}) we 
infer that we can take $\alpha_{+}=5.12$, $\alpha_{-}=3.41$ 
and $x_0=2$ in the analogue of (\ref{zeerexpliciet}). \qed\\

\noindent {\tt Remark}. In the
proof of part a we have used the
trivial estimates
$\psi_4(x)\le x$ and $B_1(x)\le x$. Using that the integers $n$ counted by
$\psi_4(x)$ satisfy $n\equiv 1({\rm mod~}4)$ and $3\nmid n$, we obtain the
sharper estimate
\begin{equation}
\label{sharper}
\psi_4(x)\le [{x+11\over 12}]+[{x+7\over 12}]\le {x\over 6}+{7\over 6}.
\end{equation}
Similarly, some computation yields that $B_1(x)\le x/2+2$. In this way
the value 9.62 appearing in
Lemma \ref{explicit} a can still be further decreased, but we have not
carried this out. Similarly the estimates in the proof of part b can
be improved.\\

\noindent {\bf Acknowledgement}. We thank Olivier Ramar\'e for giving
us access to his preprint \cite{mertens}, some results of which were
crucial in proving our main theorem. Thanks are due to Peter Shiu
for comments regarding his paper \cite{shiu}.


\begin{thebibliography}{99}
\bibitem{berndtr} B.C. Berndt, {\it Ramanujan's notebooks}. Part IV, 
Springer-Verlag, New York, 1994.
\bibitem{berndt} B.C. Berndt and R.A. Rankin,
{\it Ramanujan: Letters and commentary}, AMS, Rhode Island, 1995.
\bibitem{bono} B.C. Berndt and K. Ono, Ramanujan's unpublished manuscript 
on the partition and tau functions with proofs and commentary, The
Andrews Festschrift (Maratea, 1998), (Eds.) D. Foata, 2001, 39-110.
\bibitem{latticebible} J.H. Conway and N.J.A. Sloane, 
{\it Sphere packings, lattices and groups}, Third edition, Grundlehren der 
Mathematischen Wissenschaften
{\bf 290}, Springer-Verlag, New York, 1999.
\bibitem{cox} D.A. Cox, {\it Primes of the
form $x\sp 2 + ny\sp 2$. Fermat, class field theory and
complex multiplication}, Wiley and Sons, Inc., New York, 1989.
\bibitem{davenport} H. Davenport, {\it Multiplicative number theory}, Third 
revised edition, Springer-Verlag, New York, 2000.
\bibitem{finch} S. Finch, Mathematical constant web pages,\\
http://www.mathcad.com/library/constants/index.htm
\bibitem{kuhnlein} S. K\"uhnlein,
Partial solution of a conjecture of Schmutz.
{\it Arch. Math. (Basel)} {\bf 67} (1996), 164-172.
\bibitem{landau} E. Landau, \"Uber die Einteilung der positiven
ganzen Zahlen in vier Klassen nach der Mindestzahl der zur ihrer
additiven Zusammensetzung erforderlichen Quadrate, {\it Archiv
der Math. und Physik} {\bf 13} (1908), 305-312.
\bibitem{moree} P. Moree, Chebyshev's bias for composite numbers with
restricted prime divisors, arXiv:math.NT/0112100, to appear in
Mathematics of Computation.
\bibitem{moreerama} P. Moree, On some claims in 
Ramanujan's `unpublished' manuscript on 
the partition and tau functions, arXiv:math.NT/0201265, submitted for 
publication.
\bibitem{cazaran} P. Moree and J. Cazaran, On a claim of
Ramanujan in his first
letter to Hardy, {\it Exposition.
Math.} {\bf 17} (1999), 289-311.
\bibitem{saradha} M.R. Murty and N. Saradha, An asymptotic formula by a method 
of Selberg,
{\it C. R. Math. Rep. Acad. Sci. Canada} {\bf 15} (1993), 273-277.
\bibitem{postnikov} A.G. Postnikov, {\it Introduction to analytic
number theory}, AMS translations of mathematical monographs 68,
AMS, Providence, Rhode Island, 1988.
\bibitem{mertens} O. Ramar\'e, Sur un th\'eor\`eme de Mertens, to appear
in Manuscripta Mathematica.
\bibitem{ramare} O. Ramar\'e and R. Rumely, Primes in arithmetic
progressions, {\it Math. Comp.} {\bf 65} (1996), 397-425.
\bibitem{conjecture} P. Schmutz Schaller, Geometry of Riemann surfaces based on
closed geodesics, {\it Bull. Amer. Math. Soc. (N.S.)} {\bf 35} (1998),
193-214.
\bibitem{student} P. Schmutz Schaller, Platonische K\"orper, Kugelpackungen
und hyperbolische Geometrie, {\it Math. Semesterber.} {\bf 47} (2000),
75-87.
\bibitem{selberg} A. Selberg, {\it Collected papers}, Vol. II, Springer-Verlag, 
Berlin, 1991.
\bibitem{serre} J.-P. Serre, Divisibilit\'e de certaines fonctions
arithm\'etiques, {\it Enseignement Math.} {\bf 22} (1976), 227-260.
\bibitem{shanks} D. Shanks, The second-order term in the
asymptotic expansion of
$B(x)$, {\it Math. Comp.} {\bf 18} (1964), 75-86.
\bibitem{shanksschmid} D. Shanks and L.P. Schmid, Variations on a theorem of
Landau, I.
{\it Math. Comp.} {\bf 20} (1966), 551-569.
\bibitem{shiu} P. Shiu, Counting sums of two squares: the Meissel-Lehmer
method, {\it Math. Comp.} {\bf 47} (1986), 351-360.
\bibitem{shiuprivate} P. Shiu, Private communication, February 2002.
\bibitem{song} J.M. Song, Sums over nonnegative multiplicative functions
over integers without large prime factors. I, {\it Acta Arith.}
{\bf 97} (2001), 329-351.
\bibitem{tenenbaum} G. Tenenbaum, 
{\it Introduction to analytic and probabilistic number theory}, 
Cambridge University Press, Cambridge, 1995.



\end{thebibliography}
\end{document}